\documentclass[12pt,twoside]{extarticle}

\usepackage[utf8]{inputenc}
\usepackage{lipsum}    
\usepackage[utf8]{inputenc}
\usepackage[english]{babel}
\usepackage{geometry}
\usepackage{fancyhdr}
\usepackage{titlesec}
\usepackage{amssymb}\usepackage{amsthm}\usepackage{amsmath}\usepackage{mathrsfs}%
\usepackage[]{enumitem}
\usepackage{comment}
\usepackage[dvipsnames]{xcolor}
\usepackage{microtype}
\usepackage{hyperref}
\usepackage{mathabx}
\renewcommand{\star}{\Asterisk}
\usepackage{mlmodern} 
\usepackage{moresize}
\usepackage{footmisc}
\usepackage{appendix}

\setenumerate{itemsep=0pt,topsep=2pt,parsep=1pt}
\setitemize{itemsep = 0pt, topsep = 0pt, parsep = 0pt}
\hypersetup{colorlinks,linkcolor={red},citecolor={blue}}

\renewcommand{\thesection}{\arabic{section}}

\numberwithin{equation}{section}

\newcommand{\uc}{\MakeLowercase}

\newcommand{\bbsizz}{\fontsize{10}{0}\selectfont\ignorespaces}
\newcommand{\siz}{\fontsize{13.5}{\baselineskip}\selectfont\ignorespaces}
\newcommand{\sizz}{\fontsize{14.5}{\baselineskip}\selectfont}

\titleformat{\chapter}
    {\Large \centering \sc \vspace*{0cm}}
    {Chapter \thechapter:}
    {5pt}
    {}
    [\vspace{1cm} \titlerule]
\titleformat{\section}
    {\large \centering \bfseries} 
    {\thesection.} 
    {5pt} 
    {} 
\titleformat{\subsection}
    [runin]
    {\bfseries}
    {\thesubsection}
    {5pt}
    {}

\renewenvironment{proof}{{\itshape Proof.}}{\hfill\qedsymbol\vspace{8pt}}

\newtheoremstyle{defbfnote}
    {8pt} 
    {8pt} 
    {} 
    {12pt} 
    {\scshape} 
    {.} 
    {.5 em} 
    {\thmname{\ignorespaces #1}\thmnumber{ #2}\thmnote{ \uc{(\ignorespaces\siz{#3})}}}
\newtheoremstyle{bfnote}
    {8pt} 
    {8pt} 
    {\itshape} 
    {12pt} 
    {\scshape} 
    {.} 
    {.5 em} 
    {\thmname{\ignorespaces #1}\thmnumber{ #2}\thmnote{ \uc{(\ignorespaces\siz{#3})}}}

\theoremstyle{bfnote}
\newtheorem{theorem}{\uc{\siz{Theorem}}}[section]
\newtheorem{corollary}[theorem]{\uc{\siz{Corollary}}}
\newtheorem{lemma}[theorem]{\uc{\siz Lemma}}
\newtheorem{proposition}[theorem]{\uc{\siz Proposition}}

\newtheorem{problem}{\uc{\siz Problem}}

\theoremstyle{defbfnote}
\newtheorem{definition}[theorem]{\uc{\siz{Definition}}}
\newtheorem{remark}[theorem]{\uc{\siz Remark}}
\newtheorem{example}[theorem]{\uc{\siz Example}}

\setlist[enumerate,1]{label={(\roman*)}}
\setlist[itemize,1]{label = {$-$}}

\renewcommand{\mathbb}{\mathbf}
\newcommand{\id}{\text{\normalfont id}}
\renewcommand{\P}{\mathcal{P}}
\newcommand{\B}{\mathcal{B}}
\newcommand{\R}{\mathbb{R}}
\newcommand{\N}{\mathbb{N}}
\newcommand{\Rn}{{{\mathbb{R}}^n}}
\newcommand{\C}{{{\mathscr{C}}}}
\newcommand{\Sp}{{{\mathscr{S}}}}
\renewcommand{\L}{{\mathscr{L}}}
\newcommand{\Ln}{{\mathscr{L}^n}}
\newcommand{\norm}[1]{\lVert#1\rVert}
\newcommand{\spt}{\text{\normalfont spt}}
\renewcommand{\d}{{\mathfrak{d}}}

\renewcommand{\div}{\text{\normalfont div}}

\renewcommand{\ae}{\text{a.e.}}
\renewcommand{\phi}{\varphi}

\newcommand{\Isop}{\text{\normalfont Isop}}
\newcommand{\tv}{\textsc{tv}}

\newcommand{\dom}{\text{\normalfont dom}}
\newcommand{\minus}{\backslash}

\newcommand{\otau}{{\overline{\tau}}}

\renewcommand{\epsilon}{\varepsilon}

\newcommand{\GMM}{\text{\normalfont GMM}}
\newcommand{\ac}{\text{\normalfont AC}}
\newcommand{\bv}{\text{\normalfont BV}}

\newcommand{\Lip}{\text{\normalfont Lip}}
\newcommand{\s}{{\P^p_\infty(\Rn)}}
\newcommand{\w}{W_\infty}

\newcommand{\csg}{\text{\normalfont CSG}}
\newcommand{\dRn}{\text{\normalfont d}_\Rn}
\newcommand{\dqp}{{\mathfrak{d}_q^p}}
\newcommand{\loc}{\text{\normalfont loc}}
\newcommand{\pl}{\text{\normalfont PL}}

\newcommand{\Crit}{\text{\normalfont Crit}}
\newcommand{\curlybrace}[1]{\left\{#1\right\}}
\renewcommand{\div}{\text{\normalfont div}}
\newcommand{\Lcorner}{\mathrel{\makebox[7pt][c]{\rule{0.4pt}{6.75pt}\rule{5.5pt}{.4pt}}}}
\renewcommand{\llcorner}{{\Lcorner}}
\renewcommand{\Tilde}{\widetilde}
\renewcommand{\hat}{\widehat}

\makeatletter
\def\thanks#1{\protected@xdef\@thanks{\@thanks
        \protect\footnotetext{#1}}}
\makeatother

\title{\uppercase{\bfseries \large Isoperimetric minimizing movements and AC curves in spaces of measures}}
\author{{\sc pietro aldrigo\textsuperscript{\textdagger}} \thanks{\textsuperscript{\textdagger}{\sc Universit\"at Bern, Mathematisches Institut (MAI), Sidlerstrasse 12, 3012 Bern, Schweiz}. Email address: \texttt{pietro.aldrigo@unibe.ch}}
\thanks{The author is supported by the Swiss National Science Foundation, grant number 200021-228012. }
}
\date{ }

\begin{document}

\maketitle

\begin{abstract}
	We define a complete metric structure on the family $\text{PL}_q^p(\mathbb{R}^n)$ of probability measures with densities in $L^p(\mathbb{R}^n)$ and finite $q$-moments. We establish the existence of generalized minimizing movements for the isoperimetric ratio and characterize absolutely continuous curves in this space through weak solutions of the continuity equation with velocity fields satisfying a first-order integral condition. We also characterize absolutely continuous curves in the $\infty$-Wasserstein space and prove a Benamou--Brenier formula for $W_\infty$.
\end{abstract}

\section[introduction]{Introduction}

The minimizing movement scheme, introduced by De Giorgi in \cite{DeGiorgi1993} and further developed by Ambrosio in \cite{Ambrosio1995}, is one of the basic variational methods for constructing evolutions from time-discrete minimization problems. Its metric formulation has become a central tool in the theory of gradient flows, especially after the development of the general theory in metric spaces and in Wasserstein spaces; we refer to \cite{ambrosio2005gradient,figalli2021invitation,santambrogio20151,Santambrogio2017} for the general background. A fundamental example is the Jordan--Kinderlehrer--Otto scheme for the Fokker--Planck equation \cite{JordanKinderlehrerOtto1998}, where the evolution is obtained by iteratively minimizing an energy penalized by the squared $W_2$-distance from the previous time step (see also \cite{CailletSantambrogio2024}).

The classical Wasserstein theory is particularly effective for displacement-convex energies, such as internal, potential, and interaction energies (see e.g. Chapter 9 of \cite{ambrosio2005gradient}). Several natural functionals in the calculus of variations, however, have a different form. This is the case for Sobolev-type and isoperimetric ratios, where a first-order quantity is normalized by a Lebesgue norm. Typical examples are
\begin{align}\label{I_1}
    \mathcal{S}_r(f):=
    \frac{\norm{\nabla f}_{L^r}}{\norm{f}_{L^{r^*}}}
\quad \text{and}\quad
    \mathrm{Isop}(f):=
    \frac{\norm{\delta f}}{\norm{f}_{L^{n/(n-1)}}},
\end{align}
where $r^*:=nr/(n-r)$ is the Sobolev conjugate of $1<r<n$ and $\norm{\delta f}$ is the total variation of $f$.
These are the scale-invariant quantities associated with the sharp Sobolev and isoperimetric inequalities. In particular, both functionals admit global minimizers, up to the natural invariances of the problem. The global minimizers of the Sobolev quotient $\mathcal{S}_r$ are the Aubin--Talenti functions, as shown in the classical works of Aubin and Talenti \cite{Aubin1976,Talenti1976}. On the other hand, the global minimizers of the isoperimetric quotient $\Isop$ are given by multiples of characteristic functions of balls, in accordance with the sharp isoperimetric inequality for sets of finite perimeter; see the works of De Giorgi and Fleming--Rishel \cite{DeGiorgi1958,FlemingRishel1960}.
Their connection with optimal transport is also well established: mass-transportation methods give proofs of sharp Sobolev and Gagliardo--Nirenberg inequalities \cite{CorderoErausquinNazaretVillani2004}, quantitative anisotropic isoperimetric inequalities \cite{FigalliMaggiPratelli2010}, and sharp stability results for anisotropic Sobolev-type inequalities in $BV$ \cite{FigalliMaggiPratelli2013}.

The interaction between first-order variational quantities and optimal transport has already appeared in several works. The Wasserstein gradient flow of the total variation and the corresponding TV--JKO scheme were studied by Carlier--Poon in \cite{CarlierPoon2019} and more recently by Lin--Santambrogio in \cite{LinSantambrogio2026}, and related Euler--Lagrange equations were later considered in the work by Chambolle--Duval--Machado \cite{ChambolleDuvalMachado2023}. 

The aim of this paper is to introduce a metric framework in which transportation
of mass and regularity of the density are both part of the topology. For $(p,q)\in [1,\infty]\times (1,\infty]$, we consider the class
\begin{align*}
    \pl_q^p(\Rn)
    :=
    \left\{\mu=f\Ln :\mu\in\P_q(\Rn),\,f\in L^p(\Rn)
    \right\},
\end{align*}
endowed with the metric
\begin{align*}
    \d_q^p(f\Ln,g\Ln):=
    W_q(f\Ln,g\Ln)+\norm{f-g}_{L^p(\Rn)}.
\end{align*}
The Wasserstein term in $\d_q^p$ controls the
displacement of mass, while the $L^p$ term controls the density in a topology
which is natural for Sobolev and isoperimetric quantities. This makes
$\pl_q^p(\mathbb R^n)$ a suitable setting for variational problems whose value
depends both on transport and on analytic properties of the density.

One of the main results of the paper is the Eulerian characterization of absolutely continuous curves in $\pl_q^p(\mathbb R^n)$ given by Theorem \ref{th_ACWqp_pNotOne}, for $1<p\leq \infty$, and by Theorem \ref{th_AC_PLOne}, for the degenerate case $p=1$. Roughly speaking, we prove that a curve $\mu_{(\cdot)}:[0,T]\to \pl_q^p(\Rn)$ is absolutely continuous with respect to $\d_q^p$ if and only if it admits a Borel velocity field $(t,x)\mapsto v_t(x)$ solving the continuity equation
\begin{align*}
    \partial_t\mu_t+\div(v_t\mu_t)=0
\end{align*}
and satisfying the first-order integrability condition
\begin{align*}
    \int_0^T
    \left(
        \norm{v_t}_{L^q(\mu_t)} +\norm{\div(v_t\mu_t)}_{L^p(\Rn)}
    \right)\,dt
    <\infty .
\end{align*}
Moreover, the metric derivative with respect to $\d_q^p$ is bounded from above by the sum $\norm{v_t}_{L^q(\mu_t)} + \norm{\div(v_t\mu_t)}_{L^p}$.

The endpoint space $(\pl^\infty_\infty(\Rn),\d_\infty^\infty)$ emerges from the study of the minimizing movements of the isoperimetric functional. Indeed, the natural domain for the isoperimetric ratio is the family $\mathcal{F}$ of bounded Borel subsets with positive measure and it is defined on $\mathcal{F}$ as
\begin{align*}
\frac{\text{Per}(\Omega)}{(\Ln(\Omega))^{(n-1)/n}},
\end{align*}
where $\text{Per}(\Omega)$ denotes the perimeter of the set $\Omega\in \mathcal{F}$. This family is embedded in the complete metric space $(\P_\infty(\Rn),W_\infty)$ via the map $\iota: \mathcal{F}\hookrightarrow \P_\infty(\Rn)$, $\iota(\Omega) := (\Ln(\Omega))^{-1}\Ln\llcorner \Omega$, thereby endowing it with a canonical notion of transportation cost, and the isoperimetric ratio of $\iota(\Omega)$ is still well-defined as $\Isop(f_\Omega)$, with $f_\Omega= (\Ln(\Omega))^{-1}\chi_\Omega$ being the density of the probability $\iota(\Omega)$ and $\Isop$ as in \eqref{I_1}. However, the subfamily $\iota(\mathcal{F})$ does not form a complete metric subspace of the $\infty$-Wasserstein space. Moreover, $\iota(\mathcal{F})$ is not convex and is not stable under the relaxed limits naturally arising in variational compactness arguments.
To resolve these issues, we convexify $\iota(\mathcal{F})$ into the set $\pl_\infty^\infty(\Rn)$ of probabilities with density in $L^\infty$ and bounded support, and incorporate the $L^\infty$-norm into the metric. The resulting metric space is complete, convex, has a canonical notion  of transportation cost and the isoperimetric ratio of its elements is well-defined. 

Another main result of the paper concerns the existence of absolutely continuous generalized minimizing movements for the isoperimetric ratio. In the  space $\pl_\infty^\infty(\mathbb R^n)$, we study the implicit scheme
\begin{align*}
    \mu_{k+1}^{\tau}\in
    \operatorname*{argmin}_{\substack{\mu\in \pl_\infty^\infty(\Rn)\\ \mu=f\Ln}}
    \left\{\Isop(f)+\frac{1}{2\tau}\left(\d_\infty^\infty(\mu,\mu_k^\tau)\right)^2\right\}, \quad \tau>0.
\end{align*}
We prove that this scheme admits absolutely continuous generalized minimizing movements (GMMs). The proof is based on the abstract theory of minimizing-movement presented in Chapter 2 of the book by Ambrosio--Gigli--Savaré \cite{ambrosio2005gradient}, together with compactness and lower semicontinuity properties provided by the metric $\d_\infty^\infty$.

The $\infty$-Wasserstein space $(\P_\infty(\Rn),W_\infty)$ plays a central role in the paper. Since $\d_\infty^\infty$ contains the uniform transportation distance $W_\infty$, the study of absolutely continuous curves in $\pl_\infty^\infty(\mathbb R^n)$ requires an analogue of the classical continuity-equation characterization of Wasserstein absolutely continuous curves. For $1<q<\infty$, this characterization is classical and follows from the dynamical theory of Wasserstein spaces; see, for instance, \cite{ ambrosio2004lecture,ambrosio2021,ambrosio2005gradient,Lisini2007,santambrogio20151}. For $q=\infty$, the corresponding statement can be regarded as the formal limit as $q\to\infty$ of this theory and is closely related to the literature on dynamical transport distances and on the $W_\infty$ distance \cite{Brasco2010,champion2008wasserstein,givens}. We give a self-contained proof of the form needed here. In particular, in the appendix we show the  Benamou--Brenier formula 
\begin{align*}
W_\infty(\mu,\nu) = \min\left\{\norm{v}_{L^\infty(\tilde{\mu})}: \begin{matrix}
	\tilde{\mu} :=\L^1\llcorner[0,1] \otimes(\mu_t)_t\hfill\\
	\partial_t\mu_t + \div(v_t\mu_t) = 0\text{ on }[0,1]\hfill\\
	\mu_{(\cdot)}\text{ narrowly continuous}\hfill\\
	\mu_0 = \mu\text{ and } \mu_1=\nu\hfill\\
	\end{matrix}
	\right\}
\end{align*}
(see \cite{BenamouBrenier2000} for the original statement with $q=2$ and \cite{ambrosio2021} for the analogous version with $1<q<\infty$),
as well as an action-minimization formula 
\begin{align*}
	W_\infty(\mu,\nu) = \min\left\{\norm{\mathscr{A}_\infty}_{L^\infty(\eta)}: \begin{matrix} \eta\in \P(\C^0([0,1];\mathscr{S})),\hfill \\ (e_0)_\#\eta = \mu,\hfill \\(e_1)_\#\eta =\nu\hfill \end{matrix}\right\},
\end{align*}
where $\mathscr{A}_\infty(\omega) = \norm{\,|\dot\omega|\,}_{L^\infty([0,1])}$ is the $\infty$-action and $(\mathscr{S},\texttt{d})$ is a general Polish geodesic space.

The paper is organized as follows. In Section \ref{sec_pre} we collect the notation and the preliminary material on measure theory, optimal transport, absolute continuity in metric and Banach spaces, and GMMs. In Section \ref{sec_PLspaces} we introduce the spaces $(\pl_q^p(\Rn),\d_p^q)$, prove their basic metric properties, and establish the existence of GMMs for the isoperimetric ratio in $\pl_\infty^\infty(\Rn)$. Section \ref{sec_WinftyAC} is devoted to absolutely continuous curves in $\P_\infty(\Rn)$. In Section \ref{sec_ACPL} we prove the main characterization theorems for absolutely continuous curves in $\pl_q^p(\Rn)$, treating separately the cases $1<p\leq\infty$ and $p=1$. Section \ref{sec_Rmk} contains final remarks and open problems. The appendix contains the dynamical formulations for $W_\infty$ as well as a superposition principle.

\vspace{3mm}

\textbf{Acknowledgments.} 
I would like to thank my supervisor, Prof. Dr. Zolt\'{a}n Balogh, for the many fruitful discussions we have had during the development of this work.

\section[notation and preliminaries]{Notation and preliminaries}\label{sec_pre}
\subsection{Measure theory.} 
\renewcommand{\d}{{\normalfont\texttt{d}}}

Let $(\Sp,\d)$ be a metric space. The Borel $\sigma$-algebra of $\Sp$ is the $\sigma$-algebra $\B(\Sp)$ generated by the open subsets of $\Sp$.
For every integer $m\geq 1$, denote by $\mathcal{M}(\Sp;\R^m)$ (resp. $\mathcal{M}(\Sp)$, if $m=1$) the family of finite $\R^m$-valued (resp. real-valued) Borel measures on $\Sp$. 
The family of Borel non-negative measures is the subset $\mathcal{M}_+(\Sp)\subseteq \mathcal{M}(\Sp)$ containing only the measures $\mu$ such that $\mu(B)\geq 0$ for every Borel subset $B\in \B(\Sp)$, and the family of Borel probability measures is the subset $\P(\Sp)\subseteq \mathcal{M}_+(\Sp)$ containing only the non-negative measures $\mu$ such that $\mu(\Sp)=1$.

The space of $\R^m$-valued measures is a vector space, and forms a Banach space when endowed with the \emph{total variation norm} $\norm{\cdot}_{\tv}$ defined by
\begin{align*}
	\norm{E}_{\tv}:=\sup\left\{\int_\Sp \langle v,dE\rangle : v\in \C^0(\Sp;\R^m),\,\norm{v}_{\C^0(\Sp;\R^m)}\leq 1\right\}.
\end{align*}
Moreover, by Riesz' representation theorem, the space $(\mathcal{M}(\Sp;\R^m),\norm{\cdot}_{\tv})$  is isometrically isomorphic to the topological dual  $(\C^0_0(\Sp;\R^m),\norm{\cdot}_{\C^0(\Sp;\R^m)})^*$, where $\C^0_0(\Sp;\R^m)$ (or simply $\C^0_0(\Sp)$, if $m=1$) denotes the space of $\R^m$-valued continuous functions that vanish at $\infty$.

We say that a sequence of probability measures $(\mu_j)_{j\geq 1}\subseteq \P(\Sp)$ converges narrowly to $\mu\in \mathcal{M}_+(\Sp)$, and write $\mu_j\xrightarrow{\text{narrow}}\mu$ if 
\begin{align*}
	\lim_{j\to\infty}\int_\Sp \psi\,d\mu_j = \int_\Sp \psi\,d\mu \quad \forall \psi\in \C^0_b(\Sp).
\end{align*}

Let $F:\Sp\to \mathscr{T}$ be a Borel map between two metric spaces $\Sp$ and $\mathscr{T}$, and let $\xi\in \mathcal{M}_+(\Sp)$. The \emph{push-forward of $\xi$ through $F$} is the measure $F_\#\xi\in \mathcal{M}_+(\mathscr{T})$ defined by
\begin{align*}
	F_\#\xi(B) := \xi(F^{-1}(B))\quad \forall B\in \B(\mathscr{T}),
\end{align*}
or equivalently 
\begin{align*}
	\int_{\mathscr{T}}\phi\,dF_\#\xi = \int_\Sp \phi\circ F\,d\xi \quad \forall \phi:\mathscr{T}\to [0,\infty] \text{ Borel measurable}.
\end{align*}

\subsection{Absolute continuity in metric and dual Banach spaces.}

\begin{definition}[Absolutely continuous curves]
	Let $(\Sp,\d)$ be a metric space, $I\subseteq \R$ be an interval and let $1\leq r\leq \infty$. A curve $\omega:I\to \Sp$ is \emph{locally $r$-absolutely continuous} (or simply \emph{absolutely continuous}, if $r=1$), and we will write $\omega\in \ac^r_\loc(I;\Sp)$ (resp. $\ac_\loc(I;\Sp)$), if there exists a function $m\in L^r_\loc(I)$ such that
	\begin{align}\label{eq_def_AC_metric}
		\d(\omega(t),\omega(s))\leq \int_s^t m(r)\,dr\quad \forall s,t\in I,\,\,  s<t.
	\end{align}
	If $I$ is compact, we write $\ac^r(I;\Sp)$ (resp. $\ac(I;\Sp)$).
\end{definition}

\begin{proposition}
	Let $(\Sp,\d)$ be a metric space, $I\subseteq \R$ be an interval and $1\leq r\leq \infty$. If $\omega\in \ac^r_\loc(I;\Sp)$, then 
	\begin{align*}
	\exists \lim_{{h\to 0}}\frac{\d(\omega(t+h),\omega(t)}{|h|}=:|\dot\omega|(t)\quad \text{for almost every }t\in I.
	\end{align*}
	Moreover, the function $|\dot\omega|$  belongs to $L^r_\loc(I)$, satisfies \eqref{eq_def_AC_metric}, and if $m\in L^r_\loc(I)$  satisfies \eqref{eq_def_AC_metric}, then $|\dot\omega|\leq m$.
\end{proposition}

Fix a Banach space $(X,\norm{\cdot}_X)$. We say that a curve $\omega:[0,T]\to X$ is \emph{differentiable almost everywhere} if there exists a function $\dot\omega:[0,T]\to X$, called \emph{derivative of $\omega$} such that 
\begin{align*}
	\lim_{{h\to 0}} \norm{\frac{\omega(t+h)-\omega(t)}{h} - \dot\omega(t)}_X=0 \quad \text{for almost every }0\leq t\leq T.
\end{align*}

Let $X$ be a separable Banach space. A function $\omega:[0,T]\to X$ is \emph{Bochner measurable} if $\omega$ is Borel measurable with respect to the topology induced by the norm $\norm{\cdot}_X$. The \emph{$r$-Bochner space} $L^r([0,T];X)$ is then the space of Bochner measurable functions $\omega:[0,T]\to X$ such that the function $t\mapsto \norm{\omega(t)}_X$ belongs to $L^r([0,T])$. 

The \emph{Bochner integral of a simple function} $\sigma:[0,T]\to X$ defined by $\sigma(t):= \sum_{k=1}^N x_k \chi_{A_k}(t)$, with $x_1,\dots,x_N\in X$ and $A_1,\dots,A_N\in \B([0,T])$ pairwise disjoint, is defined as the sum $\sum_{k=1}^N x_k\L^1(A_k)\in X$. For general curves $\omega\in L^p([0,T];X)$, the \emph{Bochner integral} is defined as the limit 
\begin{align*}
\int_0^T\omega(t)\,dt := \lim_{j\to\infty}\sum_{k=1}^{N_j} x_{k,j}\L^1(A_{k,j})\in X,
\end{align*}
where $\sigma_j:= \sum_{k=1}^{N_j} x_{k,j} \chi_{A_{k,j}}$ is any sequence of simple functions such that
\begin{align*}
\lim_{j\to\infty}\int_X\norm{\omega(t)-\sigma_j(t)}_X\,dt = 0.
\end{align*}

A Banach space $(X,\norm{\cdot}_{X})$ is called \emph{dual Banach space} if there exists a Banach space $(Y,\norm{\cdot}_Y)$ such that $(X,\norm{\cdot}_{X})$ is isometrically isomorphic to the dual $(Y,\norm{\cdot}_Y)^*$; in this case, $Y$ is called \emph{predual space of $X$}. A Banach space $X$ is called \emph{reflexive} if $X$ is isometrically isomorphic to its bi-dual $((X,\norm{\cdot}_X)^*)^*$. 

If $X$ is a dual Banach space, $Y$ is a predual of $X$ and $\langle\cdot,\cdot\rangle:X\times Y\to \R$ is the duality pairing, then the curve $\omega:[0,T]\to X$ is said to be \emph{weakly$^*$ differentiable almost everywhere} if there exists a function $\dot\omega:[0,T]\to X$, called \emph{weak$^*$ derivative of $\omega$}, such that
\begin{align*}
	\lim_{h\to 0} \langle \frac{\omega(t+h)-\omega(t)}{h}, y\rangle = \langle \dot\omega(t),y\rangle\quad \forall y\in Y\quad\text{for almost every }0\leq t\leq T.
\end{align*}

We say that $\omega:[0,T]\to X$ is \emph{weakly$^*$ measurable} if the real-valued map $t\mapsto \langle \omega(t),y\rangle$ is Borel measurable for every $y\in Y$. The \emph{weak$^*$-$L^r$ space} $L^r_{w^*}([0,T];X)$ is the space of weakly$^*$ measurable functions $\omega:[0,T]\to X$ such that the function $t\mapsto \norm{\omega(t)}_X$ belongs to $L^r([0,T])$. 
The \emph{weak$^*$ integral} of a function $\omega\in L^r_{w^*}([0,T];X)$ is an element of $X$ denoted by $\int_0^T\omega(t) dt$ with the property
\begin{align*}
\langle \int_0^T\omega(t) dt,y\rangle = \int_0^T \langle \omega(t),y\rangle \,dt\quad \forall y\in Y.
\end{align*}

We shall recall the following result (we refer to Remark 2.2 of \cite{Lisini2007}, see also \cite{Ambrosio1995,AmbrosioKirchheim2000}) 

\begin{lemma}\label{lem_Bochner}
	Let $(X,\norm{\cdot}_X)$ be a separable reflexive Banach space (resp. dual Banach space with separable predual). Then $\omega\in \ac^p([0,T];X)$ if and only if $\omega$ is differentiable (resp. weakly$^*$ differentiable) almost everywhere, its derivative (resp. weak$^*$ derivative) $\dot\omega$ belongs to the Bochner space $L^r([0,T];X)$ (resp. $L^r_{w^*}([0,T];X)$), $|\dot\omega|(t) = \norm{\dot \omega(t)}_X$ for almost every $0\leq t\leq T$ and
	\begin{align}\label{eq_Bochner_1}
	\omega(t)-\omega(s) = \int_s^t\dot\omega(r)\,dr\quad \forall 0\leq s <t\leq T,
	\end{align}
	where the integral in \eqref{eq_Bochner_1} is the Bochner's integral (resp. weak$^*$ integral).
\end{lemma}

\subsection{Minimizing movements.}
We refer to Chapter 2 of \cite{ambrosio2005gradient} for all the results presented in this subsection.

Consider the following assumptions:

\begin{enumerate}[label = (A\arabic*)]
\item \label{topAss1}	$(\Sp,\d)$ is a complete metric space;
\item \label{topAss2}	$\sigma$ is a Hausdorff topology on $\Sp$, called \emph{weak topology}, that is compatible with the metric topology, i.e. 
	\begin{enumerate}
	\item \label{topAss2a} $\sigma$ is weaker than the metric topology, and
	\item \label{topAss2b} $\d:\Sp\times\Sp\to[0,\infty)$ is $\sigma$-sequentially lower semicontinuous;
	\end{enumerate}	 
\item \label{topAss3}  $\phi:\Sp\to [0,\infty]$ is a functional such that:
\begin{enumerate}
	\item \label{topAss3a}	$\phi$ is $\sigma$-sequentially lower semicontinuous, and
	\item \label{topAss3b}	if $(x_j)_{j\geq 1}\subseteq \{\phi\leq \lambda\}$ is a $\d$-bounded sequence, then exists a subsequence $(x_{j_k})\preceq(x_j)$ and $x_*\in\mathscr{S}$ such that $x_{j_k}\xrightarrow{\sigma}x_*$.
\end{enumerate}
\end{enumerate}

\vspace{8pt}

Let $\phi:\Sp\to [0,\infty]$ be a functional as in \ref{topAss3} and define $\Phi:\Sp\times (0,\infty)\times \Sp\to [0,\infty]$
\begin{align*}
\Phi(x;\tau,\overline{x}):= \phi(x) + \frac{1}{2\tau}\d^2(x,\overline{x}).
\end{align*}
For any $\tau\ge0$, the \emph{$\tau$-resolvent operator} of $\phi$ is the set-valued functional $J_\tau[\cdot]:\Sp\to 2^{\Sp}$
	\begin{align*}
	J_\tau[x]:= \operatorname*{argmin}_\Sp \,\Phi(\cdot;\tau,x):=\left\{y \in \Sp: \Phi(y;\tau,x)= \inf_{z\in \Sp} \Phi(z;\tau,x)\right\}\subseteq \Sp.
	\end{align*}

	A \emph{partition of steps} is a family $\otau:=\{\tau_j\}_{j\geq 1}\subseteq (0,\infty)$ such that:
\begin{enumerate}
	\item $|\otau|:=\sup_{j\geq 1}\tau_j\le\infty$;
	\item $\sum_{j\geq 1}\tau_j = \infty$.
\end{enumerate}	   
Given a partition of steps $\otau$, the corresponding \emph{partition of times} is the collection $\{t_j^\otau\}_{j\geq 0}\subseteq[0,\infty)$ of the numbers defined by
\begin{align*}
	t_0^\otau :=0,\quad t^\otau_j:=\sum_{k=1}^j\tau_k \quad \forall j\geq 1.
\end{align*} 
For every $j\geq 1$, the \emph{$j$-th interval} associated with the partition of steps $\otau$ is the interval $I^\otau_j:=(t^\otau_{j-1},t^\otau_j]$.

	\begin{definition}[discrete solution]
	Let $\otau$ be a partition of steps and $\overline{x}\in \Sp$ be fixed. A \emph{discrete solution} associated with $\otau$ with initial datum $\overline{x}$ is a piecewise constant function $x^\otau_{(\cdot)}:[0,\infty)\to \Sp$ such that
	\begin{align*}
	\begin{matrix}
	x^\otau_0 = \overline{x}\hfill\\
	x^\otau_t\equiv \overline x^\otau_1 \in J_{\tau_1}[\overline{x}]\quad \forall t\in I^\otau_1,\hfill\\
	x^\otau_t\equiv \overline  x^\otau_j\in J_{\tau_j}[\overline x^\otau_{j-1}] \quad\forall t\in I^\otau_j \,\, \forall j\geq 2.
	\end{matrix}
	\end{align*}
	\end{definition}

\begin{definition}[Generalized minimizing movement]
	Let $x_0\in\Sp$ be fixed. A curve $x{(\cdot)}:[0,\infty)\to \Sp$ is a \emph{generalized minimizing movement for $\phi$ starting from $x_0$} and we write $x{(\cdot)}\in \GMM(\phi;x_0)$ if there exists a sequence of partition of steps $(\otau_k)_{k\geq 1}$ with $\lim_{k\to\infty}|\otau_k|=0$ and discrete solutions  $x^{\otau_k}_{(\cdot)}$ associated with $\otau_k$ and with initial datum $\overline{x}^{\otau_k}$ such that
	\begin{gather*}
	\begin{matrix}
	\lim_{k\to\infty} \phi(\overline{x}^{\otau_k}) = \phi(x_0),\hfill\\
 	\limsup_{k\to \infty}\d(x_0^{\otau_k},x_0)<\infty,\hfill\\ 
 	x^{\otau_k}_t\xrightarrow[k\to\infty]{\sigma}x(t) \quad \forall 0\leq t\le\infty.	
	\end{matrix}
	\end{gather*}
\end{definition}

\begin{theorem}\label{th_GenreralExistence}
	If the topological assumptions \ref{topAss1}, \ref{topAss2} and \ref{topAss3} hold, then for every $x_0\in \dom\,\phi$ there exists a generalized minimizing movement $x(\cdot)\in \GMM(\phi;x_0)\cap \ac^2_\loc([0,\infty);\mathscr{S})$ such that $x(0)=x_0$.
\end{theorem}

\subsection{Wasserstein spaces.} Let $n\geq 1$ be a fixed dimension. We denote by $(\Rn,\dRn)$ the Euclidean $n$-dimensional metric space. Consider the projections $\pi_j:\Rn\times\Rn\to \Rn$ defined by $\pi_j(x_1,x_2):=x_j$ for $j\in \{1,2\}$. Given two probability measures $\mu,\nu\in \P(\Rn)$ we define the set of \emph{couplings of $\mu$ and $\nu$} as
\begin{align*}
	\Gamma(\mu,\nu):=\left\{\gamma\in \P(\Rn\times \Rn): (\pi_1)_\#\gamma  = \mu,\, (\pi_2)_\#\gamma = \nu\right\}.
\end{align*} 

\begin{definition}[Wasserstein space]
	Let $1\leq q <\infty$. The $q$-Wasserstein space is defined as the space $(\P_q(\Rn),W_q)$, where 
	\begin{gather*}
	\P_q(\Rn):=\left\{\mu\in \P(\Rn) : \int_\Rn |x|^q\,d\mu(x)<\infty\right\},\\
	W_q(\mu,\nu):=\inf_{\gamma \in \Gamma(\mu,\nu)}\left\{\left(\int_{\Rn\times \Rn} \dRn^q(x,y)\,d\gamma(x,y)\right)^\frac{1}{q}\right\}= \inf_{\gamma\in \Gamma(\mu,\nu)}\norm{\dRn(\cdot,\cdot)}_{L^q(\gamma)}.
	\end{gather*}
	The $\infty$-Wasserstein space is the space $(\P_\infty(\Rn),W_\infty)$, where
	\begin{gather*}
	\P_\infty(\Rn):=\curlybrace{\mu\in \P(\Rn): \spt\, \mu \text{ is  bounded}},\\ W_\infty(\mu,\nu):=\inf_{\gamma\in \Gamma(\mu,\nu)}\norm{\dRn(\cdot,\cdot)}_{L^\infty(\gamma)} .
	\end{gather*}
\end{definition}

It turns out (\cite{ambrosio2021,ambrosio2005gradient,figalli2021invitation}, and \cite{champion2008wasserstein,givens} for the case $q=\infty$) that for every $1\leq q\leq \infty$ $(\P_q(\Rn),W_q)$ is a complete metric space and the infimum in $W_q(\mu,\nu)$ is attained for every $\mu,\nu\in \P_q(\Rn)$. In the sequel, we will denote by $\Gamma_q(\mu,\nu)$ the set of couplings in $\gamma\in \Gamma(\mu,\nu)$ such that
\begin{align*}
	W_q(\mu,\nu) =\norm{\dRn(\cdot,\cdot)}_{L^q(\gamma)},
\end{align*}
for every $1\leq q \leq\infty$. The elements of $\Gamma_q(\mu,\nu)$ will be also called \emph{$q$-optimal couplings of $\mu$ and $\nu$}.

The topology induced by $W_q$ is in general stronger than the one of the narrow convergence of measures, namely $W_q$-convergence always implies narrow convergence of measure. If $1\leq q <\infty$ and $(\mu_j)_{j\geq 1}\subseteq \P_q(\Rn)$ is equi-compactly supported (i.e. there exists $K\subseteq \Rn$ compact such that $\spt\,\mu_j\subseteq K$ for all $j\geq 1$), then narrow convergence and $W_q$-convergence coincide. However, there exist equi-compactly supported sequences $(\mu_j)_{j\geq 1}\subseteq \P_\infty(\Rn)$ that converge narrowly but don't admit $W_\infty$-limit (an easy example is the sequence $\mu_j:= (1-j^{-1})\delta_0 + j^{-1}\delta_1$ in $\R$). 

A useful characterization of the $W_\infty$-metric is the following:
\begin{align}\label{eq_Winfty}
	W_\infty(\mu,\nu)= \inf\{\epsilon\ge0 : \mu(A)\leq \nu(A_\epsilon)\,\,\forall A\in \B(\Rn)\}, 
\end{align}
where $A_\epsilon:= \{y\in \Rn: \dRn(y,A)\le\epsilon\}$ is the $\epsilon$-neighborhood of $A$.

\subsection{Absolute continuity in $\P_q(\Rn)$ for $1<q<\infty$.} 

\begin{definition}[weak solutions of the continuity equation]\label{def_WeakSol}
	A couple $(\mu_{(\cdot)},v_{(\cdot)})$ is called \emph{weak solution of the continuity equation in $[0,T]$} (or simply \emph{solution of the continuity equation}) if $\mu_{(\cdot)}:[0,T]\to \P(\Rn)$, $[0,T]\times \Rn \ni (t,x)\mapsto v_t(x)\in \Rn$ is a Borel map, the integrability condition
\begin{align}\label{def_contEqIntegrability}
\int_0^T\norm{v_t}_{L^1(\mu_t)}\,dt \le \infty
\end{align}	
holds
	 and the partial differential equation
	 \begin{align*}
	 \partial_t \mu_t + \div(v_t\mu_t) = 0
	 \end{align*}
	 is satisfied in the distributional sense, namely
	\begin{align}\label{def_contEq}
	\int_0^T\int_\Rn\,\left(\partial_t \phi(t,x) + \langle \nabla \phi(t,x),v_t(x)\rangle\right)\,d\mu_t(x)\,dt = 0 \quad \forall \phi\in \C^\infty_c((0,T)\times \Rn). 
	\end{align}
	In this case, we call $v_{(\cdot)}$ the \emph{velocity-field of $\mu_{(\cdot)}$}.
\end{definition}

\begin{remark}\label{rmk_NarrowContWeakSol}
In Lemma 8.1.2 of \cite{ambrosio2005gradient}, the authors proved that for every weak solution $(\mu_{(\cdot)},v_{(\cdot)})$ of the continuity equation in $[0,T]$, up to modifying $\mu_t$ is a $\L^1$-negligible subset of $[0,T]$, we can suppose $\mu_{(\cdot)}:[0,T]\to \P_\infty(\Rn)$ to be narrowly continuous. Moreover, if $\mu_{(\cdot)}:[0,T]\to \P(\Rn)$ is narrowly continuous, $v_t\in L^1(\mu_t)$ for almost every $0\leq t\leq T$ and $t\mapsto \norm{v_t}_{L^1(\mu_t)}$ belongs to $L^1([0,T])$,  then $(\mu_{(\cdot)},v_{(\cdot)})$ is a  weak solution of the continuity equation in $[0,T]$ if and only if $\psi\in \C^1_c([0,T]\times \Rn)$ 
\begin{align}\label{eq_NarrowContWeakSol}
	\int_\Rn \psi(t,\cdot)d\mu_t - \int_\Rn \psi(s,\cdot)\,d\mu_s = \int_s^t\int_\Rn \left(\partial_t\psi(r,\cdot) + \langle \nabla\psi(r,\cdot),v_r\rangle\right)\,d\mu_r\,dr
\end{align}
holds for every $0\leq s\leq t\leq T$, or equivalently (Proposition 16.3 in \cite{ambrosio2021}) for every function $g\in \C^1_c(\Rn)$, the map $t\mapsto \int_\Rn g\,d\mu_t$ belongs to $\ac([0,T])$ and its derivative is
\begin{align}\label{eq_NarrowContWeakSol_2}
	\frac{d}{dt}\int_\Rn g\,d\mu_t = \int_\Rn \langle \nabla g,v_t \rangle\,d\mu_t
\end{align}
for almost every $0\leq t\leq T$. 
\end{remark}

\begin{theorem}\label{th_ACinPq}
	Let $1<q\le \infty$ and let $\mu_{(\cdot)}:[0,T]\to \P_q(\Rn)$ be a narrowly continuous curve. 
	\begin{enumerate}
	\item Suppose the existence of a Borel time-dependent vector field $(t,x)\mapsto v_t(x)$ such that $(\mu_{(\cdot)},v_{(\cdot)})$ is a weak solution of the continuity equation and the function $t\mapsto \norm{v_t}_{L^q(\mu_t)}$ belongs to $L^1([0,T])$. Then $\mu_{(\cdot)}$ is $W_q$-absolutely continuous and $|\dot \mu|^{W_q}(t) \leq \norm{v_t}_{L^q(\mu_t)}
	$ for almost every $0\leq t \leq T$, where $|\dot\mu|^{W_q}$ is the $W_q$-metric derivative of $\mu$.
	\item Suppose  $\mu_{(\cdot)}\in \ac([0,T];\P_q(\Rn))$. Then there exists a Borel time-dependent vector field $(t,x)\mapsto v_t(x)$ such that $(\mu_{(\cdot)},v_{(\cdot)})$ is a weak solution of the continuity equation and $\norm{v_t}_{L^q(\mu_t)} = |\dot\mu|^{W_q}(t)$ for almost every $0\leq t\leq T$.
	\end{enumerate}
\end{theorem}

\section[the {\normalfont $\text{\sc pl}_q^p$}-spaces and applications]{The {\normalfont $\pl_q^p$}-spaces and applications}\label{sec_PLspaces}

\renewcommand{\d}{{\mathfrak{d}}}

\begin{definition}
	Let $1\leq p,q\leq \infty$. The $\pl_q^p$-space on $\Rn$ is the set
	\begin{align*}
	\pl_q^p(\Rn):=\left\{\mu\in \P_q(\Rn): \mu\ll \Ln\text{ and }\frac{d\mu}{d\Ln}\in L^p(\Rn)\right\},
\end{align*}
	endowed	with the metric
	\begin{align*}
	\d_q^p(\mu,\nu):= W_q(\mu,\nu) + \norm{\frac{d\mu}{d\Ln}-\frac{d\nu}{d\Ln}}_{L^p}.
	\end{align*}
\end{definition}

\begin{proposition}\label{prop_topAss1}
	The space $(\pl^p_q(\Rn),\d_q^p)$ is a complete metric space.
\end{proposition}

\begin{proof}
	Let $(\mu_j)_{j\geq 1}\subseteq \pl^p_q(\Rn)$ be a $\d_q^p$-Cauchy sequence. Then $(\mu_j)_{j\geq 1}\subseteq \P_q(\Rn)$ is $W_q$-Cauchy and $(d\mu_j/d\Ln)_{j\geq 1}\subseteq L^p(\Rn)$ is $L^p$-Cauchy. By completeness of $(\P_q(\Rn),W_q)$ and of $(L^p(\Rn),\norm{\cdot}_{L^p})$ respectively, there exist $\mu\in \P_q(\Rn)$ and $f\in L^p(\Rn)$ such that
	\begin{align*}
	\lim_{j\to\infty} W_q (\mu_j,\mu) = 0 \quad \text{and}\quad \lim_{j\to\infty}\norm{\frac{d\mu_j}{d\Ln}-f}_{L^p}=0.
	\end{align*}
	As $W_q$-convergence implies weak$^*$-convergence of measures, for every test function $\psi \in \C^\infty_c(\Rn)$ we have
	\begin{align*}
	\int_\Rn \psi\,d\mu = \lim_{j\to\infty}\int_\Rn \psi \,d\mu_j = \lim_{j\to \infty}\int_\Rn \psi\frac{d\mu_j}{d\Ln}\,d\Ln = \int_\Rn \psi\,f\,d\Ln.
	\end{align*}
	Therefore $\mu = f\Ln\in \pl^p_q(\Rn)$ and $\dqp(\mu_j,\mu)\to 0$ as $j\to \infty$. 
	
\end{proof}

\begin{remark}
	If $\mu=f\Ln\in \pl^p_q(\Rn)$, then $f\in L^r(\Rn)$ for every $1\leq r\leq p$ and $\mu\in \P_s(\Rn)$ for every $1\leq s\leq q$. Therefore we have the inclusion
	\begin{align*}
	\pl_q^p(\Rn)\subseteq \pl_s^r(\Rn)\quad \forall (r,s)\in [1,p]\times[1,q]. 
	\end{align*}
	On the other hand, 
	\begin{align*}
	\pl_q^p(\Rn) \not\subseteq \pl_s^r(\Rn)\quad \text{if }r>p\text{ or }s>q.
	\end{align*}
	Indeed, if $r>p\geq 1$ and $\alpha\in (n/r,n/p)$, where $n/ \infty:=0$. Then the measure $\mu_\alpha:= f_\alpha\Ln$, with $f_\alpha(x):=C_\alpha |x|^{-\alpha}\chi_{B_1}(x)$ and $C_\alpha>0$ being the normalizing constant, belongs to $\pl_\infty^p(\Rn)$, and thus in $\pl_q^p(\Rn)$, but not in $\pl^r_s(\Rn)$ for any choice of parameters $q,s\in [1,\infty]$.
	 Similarly, if $s>q\geq 1$ and $\beta \in (q+n,s+n)$, where $\infty+n:=\infty$. Then the measure $\mu_\beta:=g_\beta\Ln$, with $g_\beta(x):= C_\beta |x|^{-\beta}\chi_{\Rn\minus B_1}(x)$, belongs to $\pl_q^\infty(\Rn)$, and thus in $\pl_q^p(\Rn)$, but not in $\pl_s^r(\Rn)$ for any choice of parameters $p,r\in [1,\infty]$.
\end{remark}

\renewcommand{\s}{{\pl}_\infty^\infty(\Rn)}
\renewcommand{\d}{\mathfrak{d}_\infty^\infty}
We now present the main motivation for the development of the theory of $\pl_q^p$ space.
Consider the isoperimetric functional $\Isop:\pl_\infty^\infty(\Rn)\to [0,\infty]$ defined by
\begin{align}\label{def_Isop}
\Isop(f\Ln):= \frac{\norm{\delta f}}{\norm{f}_{L^{n/(n-1)}}},
\end{align}
where
\begin{align*}
\norm{\delta f}:=\sup\left\{\int_\Rn f\div\,\Psi\,d\Ln: \Psi\in \C^1_c(\Rn;\Rn),\,\norm{\Psi}_{\C^0}\leq 1\right\}
\end{align*}
is the total variation of $f$. Observe that for every $f\Ln\in \s$, $f\in L^p(\Rn)$ for every $1\leq p\leq \infty$. In particular $f\in L^{n/(n-1)}(\Rn)$ and therefore the denominator of $\Isop(f\Ln)$ is always finite.

 Finally, the weak topology $\sigma$ that we consider is the topology of $L^{\frac{n}{n-1}}$-convergence of densities, namely the topology defined by the following condition
\begin{align*}
\mu_j\xrightarrow{\sigma}\mu :\iff \lim_{j\to\infty}\norm{\frac{d\mu_j}{d\Ln}-\frac{d\mu}{d\Ln}}_{L^{n/(n-1)}} =0.  
\end{align*}

\begin{proposition}\label{prop_topAss2}
	The topology $\sigma$ is weaker than the metric topology of $(\s,\d)$ and $\d$ is $\sigma$-sequentially lower semicontinuous.
\end{proposition}

\begin{proof}
	If $(\mu_j)_{j\geq 1}\subseteq \s$ is a $\d$-converging sequence, and $\mu$ is its $\d$-limit, then \eqref{eq_Winfty} implies that $(\mu_j)_{j\geq 1}$ is equi-compactly supported. By standard interpolation of Lebesgue spaces, it follows that $(d\mu_j/d\Ln)_{j\geq 1}$ converges to $d\mu/d\Ln$ in $L^{\frac{n}{n-1}}$. Thus $\mu_j\xrightarrow{\sigma}\mu$. This proves the first statement.
	
	Let $\mu_j,\nu_j,\mu,\nu\in \s$ for all $j\geq 1$ and suppose $\mu_j\xrightarrow{\sigma}\mu$ and $\nu_j\xrightarrow{\sigma}\nu$. We need to show
	\begin{align}\label{eq_assTop2_1}
	\d(\mu,\nu)\leq \liminf_{j\to \infty} \d(\mu_j,\nu_j).
	\end{align}
	Without loss of generality, up to the extraction of a subsequence we may suppose the limit to be attained and finite, and, up to the extraction of a further subsequence, we may also suppose that
	\begin{align*}
	\left(\frac{d\mu_j}{d\Ln}(x),\frac{d\nu_j}{d\Ln}(x)\right) \to \left(\frac{d\mu}{d\Ln}(x),\frac{d\nu}{d\Ln}(x)\right)\,\, \ae\,x\in \Rn.
	\end{align*}
	By Fatou's lemma for the $L^\infty$ norm, we deduce
	\begin{align}\label{eq_assTop2_2}
	\norm{\frac{d\mu}{d\Ln}-\frac{d\nu}{d\Ln}}_{L^\infty}  \leq \liminf_{j\to\infty} \norm{\frac{d\mu_j}{d\Ln}-\frac{d\nu_j}{d\Ln}}_{L^\infty}.
	\end{align}
	Define
	\begin{align*}
	\lambda := \liminf_{j\to \infty} \w(\mu_j,\nu_j).
	\end{align*}
	Then, by \eqref{eq_Winfty}, for every $\epsilon \ge0$ there exists a subsequence $({j_k})_{k\geq 1}$ such that
	\begin{align*}
	\mu_{j_k}(A)\leq \nu_{j_k}(A_{\lambda + \epsilon}) 
	\quad \forall A\in \B(\Rn)\quad \forall k\geq 1.
	\end{align*}
	By absolute continuity and $L^\frac{n}{n-1}$-convergence of the densities, it follows that 
	\begin{align}\label{eq:boh}
	\mu(B) = \lim_{{k}\to \infty}\mu_{j_k}(B) \leq \lim_{k\to\infty} \nu_{j_k}(B_{\lambda + \epsilon}) = \nu(B_{\lambda+\epsilon})
	\end{align}
	for every bounded Borel set $B\in \B(\Rn)$. Because $\spt\,\mu \cup \spt\,\nu$ is bounded, inequality \eqref{eq:boh} holds for every Borel set $B\in \B(\Rn)$. This in turn implies $\w(\mu,\nu) \leq \lambda + \epsilon$ by \eqref{eq_Winfty}. Hence, by arbitrariness of $\epsilon\ge0$, we deduce
	\begin{align}\label{eq_assTop2_3}
	\w(\mu,\nu) \leq \liminf_{j\to\infty} \w(\mu_j,\nu_j).
	\end{align}
	The claim \eqref{eq_assTop2_1} immediately follows from \eqref{eq_assTop2_2} and \eqref{eq_assTop2_3}. 
	
\end{proof}

\begin{proposition}\label{prop_topAss3}
	The functional $\Isop$ defined in \eqref{def_Isop} is $\sigma$-sequentially lower semicontinuous and, for every $\lambda\ge0$ and $\d$-bounded sequence $(\mu_j)_{j\geq 1}\subseteq \{\Isop\leq \lambda \}$ there exists $\mu_*\in \s$ and a subsequence $(\mu_{j_k})_{k\geq 1}\preceq (\mu_j)_{j\geq 1}$ such that $\mu_{j_k}\xrightarrow{\sigma}\mu_*$.
\end{proposition}

\begin{proof}
	The $\sigma$-lower semicontinuity of $\Isop$ is a direct consequence of the lower semicontinuity of the total variation under $L^1$-convergence \cite{EvansGariepy2015}. Therefore, only the second claim requires a proof.
	
	Let $\lambda \ge0$ and let $(\mu_j)_{j\geq 1}\subseteq\{\Isop\leq \lambda\}$ be a $\d$-bounded sequence. Let $K\subseteq \Rn$ be a regular compact subset such that $\spt\,\mu_j\subseteq K$ for every $j\geq 1$. Define $f_j := d\mu_j/d\Ln$ and 
	\begin{align*}
	u_j := \frac{f_j}{\norm{f_j}_{L^{n/(n-1)}}}\quad \forall j\geq 1.
	\end{align*}
	Then $\spt\,u_j\subseteq K$ and, by H\"older's inequality
	\begin{align*}
	\norm{u_j}_{\bv} = \int_\Rn u_j\,d\Ln + \norm{\delta u_j} \leq \left(\Ln(K)\right)^\frac{1}{n} \norm{u_j}_{L^{n/(n-1)}}+ \Isop(\mu_j) \leq \left(\Ln(K)\right)^\frac{1}{n} + \lambda.
	\end{align*}
	Therefore, $(u_j)_{j\geq 1}\subseteq \bv(K)$ is a bounded sequence. By virtue of the compact embedding $\bv(K)\hookrightarrow L^1(K)$ (Theorem 4 in Section 5 of  \cite{EvansGariepy2015}), there exists a function $u_*\in \bv(K)$ and a subsequence $(u_{j_k})_{k\geq 1}\preceq (u_j)_{j\geq 1}$ such that 
	\begin{align*}
	\lim_{k\to \infty}\norm{u_{j_k}-u_*}_{L^1}=0.
	\end{align*}
	Up to the extraction of a further subsequence and the extension of $u_{j_k}$ and $u_*$ by zero in $\Rn\minus K$, we can also suppose
	\begin{align}\label{eq_topAss3_1}
	\lim_{k\to\infty} u_{j_k}(x) = u_*(x)\quad \ae\,x\in \Rn.
\end{align}		
	Observe that $\d$-boundedness of $(\mu_j)_{j\geq 1}$ implies boundedness of the sequence $(\norm{f_{j_k}}_{L^{\infty}})_{k\geq 1}\subseteq (0,\infty)$, namely the existence of $M\ge0$ such that
	\begin{align}\label{eq_topAss3_3}
	\sup_{k\geq 1}\norm{f_{j_k}}_{L^{\infty}}\leq M.
	\end{align}
	 This in turn gives a finite upper bound for the sequence $(\norm{f_{j_k}}_{L^{n/(n-1)}})_{k\geq 1}\subseteq (0,\infty)$. On the other hand, as $\mu_{j_k}$ is a probability measure that is supported on the compact set $K$, H\"older's inequality yields
	 \begin{align*}
	 \norm{f_{j_k}}_{L^{n/(n-1)}}\geq \left(\Ln(K)\right)^{-\frac{1}{n}}
	 \end{align*}
	for every $k\geq 1$. Therefore, it is not restrictive to suppose 
	\begin{align}\label{eq_topAss3_2}
	\exists L:=\lim_{k\to\infty}\norm{f_{j_k}}_{L^{n/(n-1)}}\in (0,\infty).
	\end{align}
	Define $\mu^*:=f^*\Ln$, where $f_*:=Lu_*$. Using \eqref{eq_topAss3_1} and \eqref{eq_topAss3_2} we deduce $ f_{j_k}(x) \to f_*(x)$ as $k\to\infty$ for almost every $x\in \Rn$.
	Therefore, by \eqref{eq_topAss3_3} and Fatou's lemma for $L^\infty$, it follows that \begin{align}\label{eq_LpBound2}
	\norm{f_*}_{L^\infty}\leq M.
	\end{align}
	Moreover,
	\begin{align*}
	\int_\Rn f^*\,d\Ln = L \int_\Rn u^*\,d\Ln = L\,\lim_{k\to\infty}\frac{\int_\Rn f_{j_k}\,d\Ln}{\norm{f_{j_k}}_{L^{n/(n-1)}}} = 1,
\end{align*}		
	thus $\mu_*\in \pl_\infty^\infty(\Rn)$

	Fix an arbitrary $\epsilon\ge0$.
	By Egoroff's theorem, there exists a measurable subset $\Omega\subseteq K$ with $\Ln(K\minus \Omega)\le\epsilon$ such that $(f_{j_k})_{k\geq 1}$ converges uniformly to $f_*$ in $\Omega$. Let $\delta\ge0$ be fixed and let $k^*:= k^*(\Omega,\delta)\in \N$ be such that
	\begin{align*}
	\sup_{k\geq k^*}\sup_{\Omega} |f_{j_k}-f_*|\leq \delta.
	\end{align*}
	Then, by H\"older's inequality and the $L^\infty$-bounds \eqref{eq_topAss3_3} and \eqref{eq_LpBound2}, we deduce
	\begin{align}\label{eq_breaks}
	\begin{split}
	\int_{\Rn} |f_{j_k}-f_*|^\frac{n}{n-1}\,d\Ln 
	&\leq \int_{\Omega}|f_{j_k}-f_*|^\frac{n}{n-1}\,d\Ln + \int_{K\minus \Omega}\left(|f_{j_k}|+|f^*|\right)^\frac{n}{n-1}\,d\Ln\\
	&\leq \delta^\frac{n}{n-1}\Ln(K) + (2M)^\frac{n}{n-1}\epsilon
\end{split}	
	\end{align}
	for every $k\geq k^*$. Passing to the limit as $k\to \infty$ and recalling the arbitrariness of $\epsilon$, we obtain
	\begin{align*}
	\lim_{k\to\infty} \norm{f_{j_k}-f_*}_{L^{n/(n-1)}} = 0.
	\end{align*}
	Therefore, $\mu_{j_k}\xrightarrow{\sigma}\mu_*$.
	
\end{proof}

\begin{corollary}\label{cor_existenceGMM}
	For every non-negative compactly supported function $g\in L^\infty(\Rn)$ with finite total variation, there exists a generalized minimizing movement $\mu(\cdot)\in \GMM(\Isop;\mu_g)\cap \ac^2_\loc([0,\infty);\s)$, where $\mu_g:=g/\norm{g}_{L^1}\Ln$. In particular, for every bounded set $\Omega\subseteq\Rn$ of positive measure and finite perimeter, there exists a generalized minimizing movement $\mu(\cdot)\in \GMM(\Isop;\mu_{\Omega})\cap \ac^2_\loc([0,\infty);\s)$, where $\mu_\Omega := (\Ln(\Omega))^{-1}\Ln\llcorner \Omega$.
\end{corollary}

\begin{proof}
	Let $g\in L^\infty(\Rn)$ be a non-negative compactly supported function with finite total variation. Define $f:=g/\norm{g}_{L^1}$, so that $\mu_g = f\Ln\in \dom\,\Isop\subseteq \s$. Observe that Proposition \ref{prop_topAss1}, Proposition \ref{prop_topAss2} and  Proposition \ref{prop_topAss3} guarantee respectively the topological assumptions \ref{topAss1}, \ref{topAss2} and \ref{topAss3} for the metric space $(\s,\d)$ endowed with the weak topology $\sigma$ of the $L^\frac{n}{n-1}$-convergence of densities and the functional $\Isop:\s\to[0,\infty]$ defined in \eqref{def_Isop}. Therefore, Theorem \ref{th_GenreralExistence} applies and a generalized minimizing movement $\mu(\cdot)\in\GMM(\Isop;\mu_g)\cap \ac^2_\loc([0,\infty);\pl_\infty^\infty(\Rn))$ exists. 
	
\end{proof}

\begin{remark}
\newcommand{\n}{{n/(n-1)}}
\newcommand{\ninv}{{(n-1)/n}}
	The definition \eqref{def_Isop} of isoperimetric  functional $\Isop$ can actually be extended to $\pl^p_\infty(\Rn)$ for every $p\geq n/(n-1)$. It is an easy exercise to check that Proposition \ref{prop_topAss2}, Proposition \ref{prop_topAss3} and therefore Corollary \ref{cor_existenceGMM} hold true if the metric space $(\s,\d)$ is replaced by $(\pl_\infty^p(\Rn),\mathfrak{d}_\infty^p)$ for every $n/(n-1)<p<\infty$. 
	
	On the other hand, Proposition \ref{prop_topAss3} breaks down for the limit case $p=n/(n-1)$. More precisely, the only key step in the proof that fails is the inequality  \eqref{eq_breaks}. Indeed, while it is true that convergence almost everywhere combined with boundedness in $L^p$ implies, by Egoroff's theorem, strong convergence in $L^r$ for any $1\leq r<p$, it is not in general true that strong convergence in $L^p$ holds. 
	In fact it is not in general true that $\mathfrak{d}_\infty^\n$-bounded sequences $(\mu_j)_{j\geq 1}\subseteq \pl_\infty^\n(\Rn)$ contained in sublevels of $\Isop$ are precompact in the topology $\sigma$ of the strong $L^\n$-convergence. A counter-example is the following. 
	
	Let $K:=[-1/2,1/2]^n\subseteq \Rn$ and $\psi\in \C^\infty_c(B_{1/2})$ be a non-negative function with $\int_\Rn \psi\,d\Ln=1$. Consider the functions $\phi_j(x):=j^{n-1}\psi(jx)$, $j\geq 1$. Then $\phi_j$ is smooth, $\int_\Rn\phi_j\,d\Ln = j^{-1}$, and $\phi_j$ is compactly supported in $K$. Define now the sequence $(\mu_j)_{j\geq 1}$ as
	\begin{align*}
	\mu_j:=f_j\Ln,\quad f_j(x):= \left(1-\frac{1}{j}\right)\chi_{K} + \phi_j(x). 
	\end{align*}
	It is indeed elementary to check that $(\mu_j)_{j\geq 1}$ is a $\mathfrak{d}^{n/(n-1)}_\infty$-bounded sequence of probability measure, with $(\norm{\delta f_j})_{j\geq 1}$ bounded but with no $\sigma$-converging subsequences. 

\end{remark}

Instead of the isoperimetric functional, one could also consider the $r$-Sobolev ratio functional $\mathcal{S}_r:\pl_\infty^{p}(\Rn)\to [0,\infty]$ defined by 
\begin{align*}
	\mathcal{S}_r(f\Ln):=\frac{\norm{\nabla f}_{L^r}}{\norm{f}_{L^{r^*}}},
\end{align*}
where $p>r^*$, and prove existence of GMMs for $\mathcal{S}_r$ in $(\pl^p_\infty(\Rn),\mathfrak{d}_\infty^p)$ by adapting the previous proofs.

\section[ac curves in the $\infty$-wasserstein space]{AC curves in the $\infty$-Wasserstein space}\label{sec_WinftyAC}

We recall and prove  the characterization of absolutely continuous curves with respect to $W_\infty$. Namely, a narrowly continuous curve $\mu_{(\cdot)}:[0,T]\to \P_\infty(\Rn)$ is absolutely continuous with respect to $W_\infty$ if and only if there exists a Borel vector field $(t,x)\mapsto v_t(x)$ such that $(\mu_{(\cdot)},v_{(\cdot)})$ solves the continuity equation in the sense of Definition \ref{def_WeakSol}, and
\begin{align*}
\int_0^T\norm{v_t}_{L^\infty(\mu_t)}\,dt<\infty.
\end{align*}
Moreover, the metric derivative is exactly the minimal admissible $L^\infty(\mu_t)$-norm of such velocity fields.

For $1<q<\infty$, the corresponding statement is classical and follows from the standard dynamic characterization of absolutely continuous curves in Wasserstein spaces \cite{ambrosio2005gradient, Lisini2007,  santambrogio20151}. The endpoint case $q=\infty$ is well known at a formal level and is often viewed as the limiting case $q\to\infty$; it is also explicitly mentioned in the literature on dynamical transport distances \cite{Brasco2010,santambrogio20151}. However, the passage from finite $q$ to $q=\infty$ is not completely trivial. For this reason, and since the endpoint characterization is used in the sequel, we include a self-contained proof. To the best of our knowledge, the precise form needed here is not available in the literature as a standalone statement with proof.

Let us start by proving some basic properties of the $\infty$-Wasserstein space.

\begin{lemma}\label{Lemma_Lsc}
	Let $(\mu_j)_{j\geq 1}\subseteq \P(\Rn)$ be a sequence of probability measures that narrowly converges to $\mu\in \P(\Rn)$. Then
	\begin{align*}
	\norm{\psi}_{L^\infty(\mu)}\leq \liminf_{j\to\infty}\norm{\psi}_{L^\infty(\mu_j)}
	\end{align*}
	for every lower semicontinuous function $\psi:\Rn\to [0,\infty]$.  
\end{lemma}

\begin{proof}
	Recall that a function $\psi:\Rn\to \R\cup\{\infty\}$ is lower semicontinuous if and only if the sets $\{\psi\ge \lambda\}\subseteq \Rn$ are open for every $\lambda \in \R$, and that  whenever $\mu_j$ narrowly converges to $\mu$, then
	\begin{align*}
	\mu(U)\leq \liminf_{j\to\infty}\mu_j(U)\quad \forall U\subseteq \Rn\text{ open}.
	\end{align*}
	
	Let $\psi:\Rn\to \R\cup\{\infty\}$ be lower semicontinuous and suppose that $\mu_j$ narrowly converges to $\mu$. By contradiction, suppose the existence of a subsequence $(\mu_{j_k})_{k\geq 1}$ such that 
	\begin{align*}
	\exists \lim_{k\to\infty}\norm{\psi}_{L^\infty(\mu_{j_k})} \le \norm{\psi}_{L^\infty(\mu)} - \delta,
	\end{align*}
	for some $\delta \ge 0$. Then there exists $k^*=k^*(\delta)\in \N$ such that
	\begin{align*}
	\mu_{j_k}\left(\curlybrace{\psi \ge \norm{\psi}_{L^\infty(\mu)}-\frac{\delta}{2}}\right) = 0 \quad \forall k\geq k^*.
	\end{align*}
	Therefore,
	\begin{align*}
	\begin{split}
	0 
	\le 
	\mu\left(\curlybrace{\psi \ge \norm{\psi}_{L^\infty(\mu)}-\frac{\delta}{2}}\right)
	\leq 
	\liminf_{j\to\infty}\mu_{j_k}\left(\curlybrace{\psi \ge \norm{\psi}_{L^\infty(\mu)}-\frac{\delta}{2}}\right) 
	=
	 0.
	\end{split}
	\end{align*}
\end{proof}

\begin{lemma}\label{lem_WinftyIsLSC}
	Let $(\mu_j)_{j\geq 1},(\nu_j)_{j\geq 1}\subseteq \P_\infty(\Rn)$ be sequences that converge narrowly to $\mu\in \P_\infty(\Rn)$ and $\nu\in \P_\infty(\Rn)$ respectively, and let $\gamma_j\in \Gamma(\mu_j,\nu_j)$ for every $j\geq 1$. 
	\begin{enumerate}
	\item\label{Lemma_LSC_i} There exists $\gamma\in \Gamma(\mu,\nu)$ and a subsequence $(\gamma_{j_k})_{k\geq 1}\preceq (\gamma_{j})_{j\geq 1}$ such that $\gamma_{j_k}$ converges narrowly to $\gamma$.
	\item\label{Lemma_LSC_ii} If $\gamma_j\in \Gamma_\infty(\mu_j,\nu_j)$  for every $j\geq 1$, then every narrow limit point $\gamma$ of $\{\gamma_j : j\geq 1\}$ satisfies
	\begin{align*}
	\norm{\text{\normalfont d}_{\Rn}(\cdot,\cdot)}_{L^\infty(\gamma)} \leq \liminf_{j\to\infty}W_{\infty}(\mu_j,\nu_j).
	\end{align*}
	\item \label{Lemma_LSC_iii} The $\infty$-Wasserstein distance is narrowly lower semicontinuous.
	\end{enumerate}
\end{lemma}

\begin{proof}
\ref{Lemma_LSC_i} The families $\{\mu_j:j\geq 1\},\{\nu_j:j\geq 1\}\subseteq \P(\Rn)$ are both tight. This implies that $\{\gamma_j : j\geq 1\}\subseteq \P(\Rn\times \Rn)$ is tight as well, and by Prokhorov's theorem we conclude.

\ref{Lemma_LSC_ii} The Euclidean distance $\dRn(\cdot,\cdot)$ is continuous. In particular, it is lower semicontinuous, therefore Lemma \ref{Lemma_Lsc} applies and the claim follows.

\ref{Lemma_LSC_iii} This immediately follows from \ref{Lemma_LSC_i} and \ref{Lemma_LSC_ii}.

\end{proof}

\renewcommand{\d}{\text{\normalfont d}}

If $(\mathscr{S},\mathfrak{d})$ is a complete metric space and $\omega:[0,T]\to (\mathscr{S},\d)$ is an absolutely continuous
then there exists a strictly increasing map $\tau:[0,S]\to [0,T]$ such that $\omega\circ\tau$ is absolutely continuous and the metric derivative of $\omega\circ \tau$ is 1 almost everywhere in $[0,S]$. This reparametrization will be called \emph{arc-length reparametrization of $\omega$}, and the continuous monotone extension of its inverse $\sigma:=\tau^{-1}:[0,T]\to [0,S]$ is absolutely continuous and enjoys the property
\begin{align*}
	\sigma'(t) = |\dot\omega|(t)\quad \text{for almost every }0\leq t\leq T
\end{align*} 
(see e.g. Lemma 1.1.4 of \cite{ambrosio2005gradient})

\begin{lemma}\label{lem_timeInvariance}
	Let $\sigma:[0,T]\to [0,S]$ be a non-decreasing, absolutely continuous surjective map. If a pair $(\nu_{(\cdot)},w_{(\cdot)})$ is a narrowly continuous weak solution for the continuity equation in $[0,S]$ then $(\mu_{(\cdot)},v_{(\cdot)}):=(\nu_{\sigma(\cdot)},\sigma'(\cdot)w_{\sigma(\cdot)})$ is a narrowly continuous weak solution for the continuity equation in $[0,T]$. 
\end{lemma}

\begin{proof}
	Clearly $t\mapsto \mu_t = \nu_{\sigma(t)}$ is narrowly continuous. The integrability condition \eqref{def_contEqIntegrability} is an immediate consequence of the change of variables formula for  absolutely continuous monotone functions. Therefore, only \eqref{def_contEq} requires a proof.
	
	Fix a test function $\phi\in \C^\infty_c((0,T)\times \Rn)$ and denote by $\phi_t(\cdot)$ the function $\phi(t,\cdot)\in \C^\infty_c(\Rn)$. Then, using \eqref{eq_NarrowContWeakSol}, the change of coordinates $s = \sigma(t+r)$ and recalling that there exists $\delta>0$ such that $\phi_t\equiv 0$ for $t\in [0,T]\minus[\delta,T-\delta]$, we obtain
	\begin{align*}
	\begin{split}
	\int_0^T\int_\Rn \partial_t\phi_t\,d\mu_t\,dt 
	&= \lim_{h\downarrow 0}\int_0^T\int_\Rn \frac{\phi_{t}-\phi_{t-h}}{h}\,d\mu_t\,dt\\
	&= -\lim_{h\downarrow 0}\int_0^T\frac{1}{h}\left(\int_\Rn \phi_t\,d\nu_{\sigma(t+h)} - \int_\Rn \phi_t\,d\nu_{\sigma(t)}\right)\,dt\\
	& = - \lim_{h\downarrow 0}\int_0^T\frac{1}{h}\int_{\sigma(t)}^{\sigma(t+h)}\int_\Rn \langle\nabla \phi_t, w_s\rangle\,d\nu_{s}\,ds\,dt\\
	&= - \lim_{h\downarrow 0}\int_0^T\frac{1}{h} \int_0^{h}\int_\Rn \langle \nabla \phi_t,v_{t+r}\rangle\,d\mu_{t+r}\,dr\,dt\\
	&= - \lim_{h\downarrow 0} \frac{1}{h}\int_0^h \int_0^T \int_\Rn \langle \nabla \phi_{t-r},v_{t}\rangle \,d\mu_{t}\,dt\,dr\\
	&= - \int_0^T\int_\Rn \langle \nabla \phi_t,v_t\rangle \,d\mu_t\,dt,
	\end{split}
	\end{align*}
	which is exactly \eqref{def_contEq}. The claim then follows by arbitrariness of the choice of $\phi$.

\end{proof}

We say that $\{\rho_\epsilon: \epsilon \ge 0\}$ is \emph{a family of strictly positive rapidly decreasing  mollifiers} if $\rho\in \C^\infty_0(\Rn)\cap L^1(\Rn)$ is a positive radially symmetric function with $\int_\Rn\rho \Ln =1$, $x\mapsto |x|^N\rho(x)$ belongs to $L^1(\Rn)$ for every $N\geq 0$ and $\rho_\epsilon = \epsilon^{-n}\rho(\cdot/\epsilon)$ for every $\epsilon \ge 0$. 

If $M\in \mathcal{M}(\Rn;\R^{m})$ is a $m$-valued Borel measure and $\rho\in \C^\infty_0(\Rn)\cap L^1(\Rn)$, we define the convolution of $M$ by $\rho$ as the function $M\star \rho:\Rn\to \R^{m}$ 
\begin{align*}
M\star \rho (x):=\int_{\Rn} \rho(x-y)\,dM(y)\quad \forall x\in \Rn.
\end{align*} 

\begin{lemma}\label{Lem_Approx}
   Let $\mu_{(\cdot)}:[0,T]\to \P_\infty(\Rn)$ be a narrowly continuous solution of the continuity equation associated with a Borel velocity field $(t,x)\mapsto v_t(x)$ that satisfies 
 \begin{align*}
 	\int_0^T \norm{v_t}_{L^\infty(\mu_t)}\,dt <\infty
\end{align*}  
 and let $\{\rho_\epsilon : \epsilon  \ge  0\}$ be a family of strictly positive rapidly decreasing mollifiers. Define 
\begin{align}\label{lem_def_approx}
	f_t^\epsilon := \mu_t \star \rho_\epsilon,\quad E_t^\epsilon := (v_t\mu_t)\star \rho_\epsilon,\quad v_t^\epsilon:=\frac{E_t^\epsilon}{f_t^\epsilon},\quad \mu_t^\epsilon := f_t^\epsilon\Ln
\end{align}
for all $0\leq t\leq T$ and every $\epsilon \ge 0$.
Then $(\mu_{(\cdot)}^\epsilon,v_{(\cdot)}^\epsilon)$ is a solution of the continuity equation for every $\epsilon \ge 0$ and 
\begin{align}
 \int_0^T \left(\sup_{U}\curlybrace{|v_t^\epsilon|} + \Lip_{U}(v_t^\epsilon)\right)\,dt \le \infty \quad & \forall U\in \B(\Rn)\text{ bounded set},\label{lem_approx_1}\\
\norm{v_t^\epsilon}_{L^\infty(\mu_t^\epsilon)}\leq \norm{v_t}_{L^\infty(\mu_t)}  \quad& \text{for almost every } 0\leq t\leq T \label{lem_approx_2}\\
 E_t^\epsilon\xrightarrow[\epsilon \downarrow 0]{\text{\normalfont narrow}}v_t\mu_t \quad& \text{for almost every } 0\leq t\leq T \label{lem_approx_3}\\
\mu_t^\epsilon \xrightarrow[\epsilon \downarrow 0]{\text{\normalfont narrow}}\mu_t \quad& \forall 0\leq t\leq T \label{lem_approx_4}.
\end{align}

\end{lemma}

\begin{proof} 
	It is easy to verify that, under the standing assumptions, $(\mu_{(\cdot)}^\epsilon,v_{(\cdot)}^\epsilon)$ is a weak solution of the continuity equation in the sense of Definition \ref{def_WeakSol}.

The proofs of properties \eqref{lem_approx_1}, \eqref{lem_approx_3} and \eqref{lem_approx_4} are done exactly as in Lemma 8.1.9 of \cite{ambrosio2005gradient}. Therefore, it is enough to prove \eqref{lem_approx_2}, and it is enough to observe that for almost every $0 \leq t \leq T$ we have
\begin{align*}
	|v_t^\epsilon(x)| = \frac{|E_t^\epsilon(x)|}{f_t^\epsilon(x)} \leq \frac{\int_\Rn \rho_\epsilon (x-y)|v_t(y)|\,d\mu_t(y)}{\int_\Rn \rho_\epsilon(x-y)\,d\mu_{t}(y)}\leq \norm{v_t}_{L^\infty(\mu_t)}\quad \forall x\in \Rn.
	\end{align*}

\end{proof}

\begin{theorem}\label{th_ineqA}
	Let $\mu_{(\cdot)}:[0,T]\to \P_\infty(\Rn)$ be a narrowly continuous solution of the continuity equation for a Borel velocity-field $(t,x)\mapsto v_t(x)$ such that $t\mapsto \norm{v_t}_{L^\infty(\mu_t)}$ belongs to $L^1([0,T])$. Then $\mu_{(\cdot)}\in \ac([0,T];\P_\infty(\Rn))$ and the $W_\infty$-metric derivative $t\mapsto |\dot\mu|^{W_\infty}(t)$ of $\mu_{(\cdot)}$ satisfies
	\begin{align*}
	|\dot\mu|^{W_\infty}(t)\leq \norm{v_t}_{L^\infty(\mu_t)}\quad \ae\,0\leq t\leq T.
	\end{align*}
\end{theorem}

\begin{proof}
	Thanks to Lemma \ref{Lem_Approx}, there are approximations $\{(\mu_{(\cdot)}^\epsilon,v_{(\cdot)}^\epsilon) : \epsilon \ge 0\}$ that are regular solutions of the continuity equation. Therefore, using Proposition 8.1.8 of \cite{ambrosio2005gradient}, if $(t,x)\mapsto T_t^\epsilon(x)$ is the flow-map associated with the time-dependent vector field $(t,x)\mapsto v_t^\epsilon(x)$, then $T^\epsilon$ is defined in $[0,T]\times \Rn$ and $\mu_t^\epsilon = (T_t^\epsilon)_\#\mu_0^\epsilon$ for every $0\leq t\leq T$.
	
	Fix $x\in \Rn$ and $0\leq t_1 \le t_2 \leq T$. Then, by definition of flow-map and \eqref{lem_approx_2}, 
	\begin{align*}
	|T_{t_2}^\epsilon(x)-T_{t_1}^\epsilon(x)| \leq \int_{t_1}^{t_2}|\dot T_s^\epsilon(x)|\,ds \leq \int_{t_1}^{t_2}\norm{v_s}_{L^\infty(\mu_s)}\,ds.
	\end{align*}
	Since $\gamma_{t_1,t_2}^\epsilon=(T_{t_2}^\epsilon \times T_{t_1}^\epsilon)_\#\mu_0^\epsilon$ is a coupling of $\mu_{t_1}^\epsilon$ and $\mu_{t_2}^\epsilon$, it follows that
	\begin{align}\label{eq_oneIneq_1}
	W_\infty(\mu_{t_1}^\epsilon,\mu_{t_2}^\epsilon)\leq \norm{\dRn(\cdot,\cdot)}_{L^\infty(\gamma_{t_1,t_2}^\epsilon)}\leq	\int_{t_1}^{t_2}\norm{v_s}_{L^\infty(\mu_s)}\,ds. 
	\end{align}
	Take the inferior limit of \eqref{eq_oneIneq_1} as $\epsilon \downarrow 0$ and combine the narrow convergence result \eqref{lem_approx_4} and the lower semicontinuity of $W_\infty$ under narrow convergence given by Lemma \ref{lem_WinftyIsLSC} to obtain
	\begin{align*}
	W_\infty(\mu_{t_1},\mu_{t_2})\leq \int_{t_1}^{t_2}\norm{v_s}_{L^\infty(\mu_s)}\,ds.
	\end{align*}
	This ends the proof.	
	
\end{proof}

\begin{theorem}\label{th_ineqB}
	Let $\mu_{(\cdot)}\in \ac([0,T];\P_\infty(\Rn))$. There exists a Borel time-dependent vector field $(t,x)\mapsto v_t(x)$ such that
	\begin{align}
	\partial_t \mu_t + \div(v_t\mu_t)= 0 \quad & \text{weakly in }[0,T]\times \Rn \label{lem_ineqB_cont}
	\\
	v_t\in L^\infty(\mu_t)\quad &\text{for almost every }0\leq t\leq T \label{lem_ineqB_Linfty}
	\\
	\norm{v_t}_{L^\infty(\mu_t)} \leq |\dot \mu|^{W_\infty}(t)\quad &\text{for almost every }0\leq t\leq T \label{lem_ineqB_third}
	\end{align}
\end{theorem}

\begin{proof}
	By virtue of Lemma \ref{lem_timeInvariance}, it is not restrictive to suppose $|\dot \mu|^{W_\infty}\equiv 1$ almost everywhere in $[0,T]$.
	
	Fix a test function $\psi\in \C^\infty_c(\Rn)$ and consider the function $t\mapsto \mu^\psi(t):=\int_\Rn\psi\,d\mu_t$. For any $0\leq s \le t\leq T$, let $\gamma_{s,t}\in \Gamma_\infty(\mu_s,\mu_t)$ be an optimal coupling, and 
 denote by  $H_\psi:\Rn\times\Rn \to [0,\infty)$ the function defined by
	\begin{align*}
	H_\psi(x,y):=\begin{cases}
	|\nabla \psi(x)| &,\text{ if }x=y,\\
	\text{ }\\
	\frac{|\psi(x)-\psi(y)|}{|x-y|} &,\text{ if }x\neq y
	\end{cases}.
	\end{align*}
	Then, for any fixed $0\le t\le T$ and $h \ge 0$ such that $0\le t+h \le T$, we have
	\begin{align*}
	\begin{split}
	\frac{|\mu^\psi(t+h)-\mu^\psi(t)|}{|h|} 
	&\leq \frac{1}{|h|}\int_{\Rn\times \Rn} |x-y| H_\psi(x,y)\,d\gamma_{t+h,t}(x,y)\\
	&\leq \frac{W_\infty(\mu_{t+h},\mu_t)}{|h|}\int_{\Rn\times\Rn} H_\psi(x,y)\,d\gamma_{t+h,t}(x,y).
	\end{split}
	\end{align*}
	Observe that, by Lemma \ref{lem_WinftyIsLSC}, for every $0\le t \le T$, $\gamma_{t+h,t}$ converges narrowly to $(\id\times \id)_\#\mu_t$ as $h\to 0$. Moreover, the function $H_\psi$ is upper semicontinuous. Therefore, for any $0\le t\le T$ such that the metric derivative $|\dot \mu|^{W_\infty}(t)$ exists and equals 1, we have
	\begin{align}\label{eq_ineqB_1}
		\limsup_{h\to 0} \frac{|\mu^\psi(t+h)-\mu^\psi(t)|}{|h|} \leq  \norm{\nabla \psi}_{L^1(\mu_t)}.
	\end{align}
	
\newcommand{\tmu}{\Tilde{\mu}}
	Define the positive measure in $\tmu:= \L^1\llcorner [0,T]\otimes (\mu_t)_t\in \mathcal{M}_+([0,T]\times \Rn)$. Then, 
	for every smooth test function $\phi\in \C^\infty_c((0,T)\times \Rn)$ we have
	\begin{align*}
	\begin{split}
		\int_{[0,T]\times \Rn} \partial_s \phi(s,x)\,d\tmu(s,x)
		&= \lim_{h\downarrow 0}\int_{[0,T]\times \Rn} \frac{\phi(s,x)-\phi(s-h,x)}{h}\,d\tmu(s,x)\\
		&= \lim_{h\downarrow 0} \int_{[0,T]} \frac{\mu^{\phi(s,\cdot)}(s)- \mu^{\phi(s,\cdot)}(s+h)}{h}\,ds.
	\end{split}
	\end{align*}
	Hence, recalling \eqref{eq_ineqB_1},
	\begin{align}\label{eq_ineqB_2}
	\begin{split}
	|\int_{[0,T]\times \Rn} \partial_s \phi(s,x)\,d\tmu(s,x)| 
	&\leq \int_{[0,T]}  \norm{\nabla \phi(s,\cdot)}_{L^1(\mu_s)}\,ds.
	\end{split}
	\end{align}
	Consider the linear subspace
	\begin{align*}
	V:=\curlybrace{\nabla \phi : \phi\in \C^\infty_c((0,T)\times \Rn)}\subseteq L^1(\tmu;\Rn)
	\end{align*}
	endowed with the restriction of the $L^1(\tmu)$-norm 
	and the linear functional $L:V\to \R$
	\begin{align*}
	L(\nabla \phi):=-\int_{[0,T]\times\Rn}\partial_s\phi(s,x)\,d\tmu(s,x).
	\end{align*}
	By \eqref{eq_ineqB_2}, $L$ is a bounded linear operator with $\Lip(L)\leq 1$. Therefore, by Hahn-Banach, $L$ can be extended to a functional $\tilde{L}:L^1(\tmu;\Rn)\to \R$ with $\Lip(\tilde{L}) = \Lip(L)\leq 1$. By Riesz representation theorem of the dual of $L^1$, there exists one unique element $v\in L^\infty(\tmu;\Rn)$ such that
	\begin{gather}\label{eq_ineqB_3}
	\tilde{L}(w)= \int_{[0,T]\times \Rn}\langle v , w\rangle \,d\tmu\quad \forall w\in L^1(\tmu;\Rn),\\
	\norm{v}_{L^\infty(\tmu)} = \Lip(\tilde{L})\leq 1.\label{eq_ineqB_4}
	\end{gather}
	
	Define $v_t := v(t,\cdot)$ for $\L^1$-almost every $0\leq t\leq T$. Then \eqref{eq_ineqB_3} trivially gives \eqref{lem_ineqB_cont} and \eqref{lem_ineqB_Linfty}, and \eqref{eq_ineqB_4} yields \eqref{lem_ineqB_third}. This ends the proof.
	
\end{proof}

Combining Theorem \ref{th_ACinPq}, Theorem \ref{th_ineqA} and Theorem \ref{th_ineqB}, one obtains the following characterization of absolute continuity in $(\P_q(\Rn),W_q)$ for every $1<q\leq \infty$ in terms of solutions of the continuity equation, that we formulate as the following unified corollary.

\begin{corollary}\label{cor_ACinPinfty}
	Let $1<q\leq \infty$.
	For curve $\mu_{(\cdot)}:[0,T]\to \P_q(\Rn)$, the following are equivalent:
	\begin{enumerate}
	\item $\mu_{(\cdot)}\in \ac([0,T];\P_q(\Rn))$;
	\item $\mu_{(\cdot)}$ is narrowly continuous, there exists a Borel vector field $(t,x)\mapsto v_t(x)$ such that $(\mu_{(\cdot)},v_{(\cdot)})$ is a solution of the continuity equation, the function $t\mapsto \norm{v_t}_{L^q(\mu_t)}$ belongs to $L^1([0,T])$ and $|\dot \mu|^{W_q}(t) = \norm{v_t}_{L^q(\mu_t)}
	$ for almost every $0\leq t \leq T$.
	\end{enumerate}
\end{corollary}

\section[absolutely continuous curves in {${\text{\sc pl}}_q^p({\text{\sc r}}^n)$}]{Absolutely continuous curves in {${\pl}_q^p(\mathbb{R}^n)$}}\label{sec_ACPL}
\newcommand{\partialf}{\mathcal{D}f}
\renewcommand{\d}{{\mathfrak{d}_\infty^p}}
\renewcommand{\s}{{\pl_q^p(\Rn)}}

Throughout this section any curve of probability measures $\mu_{(\cdot)}$ will be implicitly supposed to be defined on a fixed interval $[0,T]$, for some $T>0$, to take values in either $\mathcal{P}_q(\Rn)$ or $\s$, and to be narrowly continuous.

As in this section more than one notion of absolute continuity for a curve $\mu_{(\cdot)}$ will be considered, we introduce the following convention.  We shall say that $\mu_{(\cdot)}$ is $W_q$-absolutely continuous if $\mu_{(\cdot)}\in \ac([0,T];\P_q(\Rn))$. Analogously, we say that $\mu_{(\cdot)}$ is $\dqp$-absolutely continuous $\mu_{(\cdot)}\in \ac([0,T];\s)$. The notation $|\dot \mu|(\cdot)$ will be exclusively used for denoting the $\dqp$-metric derivative of $\mu_{(\cdot)}$, while its $W_q$-metric derivative will be denoted by $|\dot \mu|^{W_q}(\cdot)$. In a similar way, the $L^p$.metric derivative of a curve $f_{(\cdot)}:[0,T]\to L^p(\Rn)$, whenever it is defined, is denoted by $|\dot f|^{L^p}$.

\begin{remark}\label{rmk_counter}
	$\dqp$-absolute continuity is strictly stronger than $W_q$-absolute continuity if $1< p\leq \infty$. Indeed, if $\Lambda_0\in \B(\Rn)$ is a bounded subset with positive measure, $0\neq V\in \Rn$ is a fixed vector and $\Lambda_t$ is the set 
	\begin{align*}
	\Lambda_t := tV+ \Lambda_0 :=\curlybrace{x+tV \in \Rn : x\in \Lambda_0}\quad \forall 0\le t\leq T,
	\end{align*}
	then it is easy to check that the curve $t\mapsto \mu_{t}:= (\Ln(\Lambda_t))^{-1}\Ln\llcorner \Lambda_t\in \pl^\infty_\infty(\Rn)$ is $W_\infty$-absolutely continuous, and therefore $W_q$-absolutely continuous for any $1<q\leq\infty$. Indeed, $\mu_{(\cdot)}$ is narrowly continuous and a velocity-field for $\mu_{(\cdot)}$ is given by the constant vector $v_t(x)\equiv V$. Therefore, Theorem \ref{th_ineqA} guarantees the $W_\infty$-absolute continuity of $\mu_{(\cdot)}$. 
	
	On the other hand, $f_{(\cdot)}:=d\mu_{(\cdot)}/d\Ln$ fails to be $L^p$-absolutely continuous for any $1\le p\leq \infty$. More generally, if $\Omega_t\in \B(\Rn)$ is a bounded set with positive measure for every $0\leq t\leq T$ and $\mu_{(\cdot)}:[0,T]\to \s$ is the curve $\mu_{t}:=(\Ln(\Omega_t))^{-1}\Ln\llcorner \Omega_t$, then $f_{(\cdot)}:=d\mu_{(\cdot)}/d\Ln$ is $L^p$-absolutely continuous for some $1\le p\leq \infty$ if and only if $\Omega_t$ is constant for the pseudo-metric $\Ln(\cdot\triangle \cdot)$, i.e. if and only if 
	\begin{align}\label{rmk_ACLp}
	\Ln(\Omega_t\triangle \Omega_s) = 0 \quad \forall 0\leq s\leq t\leq T.
	\end{align}
	
	Clearly, if \eqref{rmk_ACLp} holds true, than $f_{(\cdot)}\equiv f_0$ is constant almost everywhere, and therefore absolutely continuous. Suppose now that the curve densities $f_{(\cdot)}$ described above is $L^p$-absolutely continuous for some $1\le p\leq \infty$. Then 
	\begin{align*}
	\Ln(\Omega_t) = 
	\begin{cases}
		{\norm{f_t}_{L^p}^{-p/(p-1)}} &, \text{ if }1<p <\infty\\
		\norm{f_t}_{L^\infty}^{-1} &,\text{ if }p=\infty
	\end{cases}\quad \text{for all }0\leq t\leq T.
	\end{align*} 
	The triangle inequality in $L^p$ and $L^p$-absolute continuity of $f_{(\cdot)}$ imply that $t\mapsto \norm{f_t}_{L^p}$ is continuous in $[0,T]$, and therefore attains both minimum and maximum. This in turn implies the existence of positive and finite constants $0\le m\le M\le \infty$ such that
	\begin{align*}
	m\le \Ln(\Omega_t)\le M\quad \forall 0\leq t\leq T.
	\end{align*}
	Fix $0\leq s\le t\leq T$. Then 
	\begin{align*}
	\norm{f_t-f_s}_{L^\infty} 
	\begin{cases}
		 \displaystyle = 0 &,\text{ if }f_t = f_s\\
		\displaystyle \geq M^{-1} &,\text{ if } f_t\neq f_s
	\end{cases},
	\end{align*}
	and therefore $f_{(\cdot)}$ is $L^\infty$-absolutely continuous if and only if it is constant in $L^\infty(\Rn)$; this is equivalent to \eqref{rmk_ACLp}.
	
	Suppose now $1\le p \le \infty$. Then
	\begin{align*}
	\norm{f_t-f_s}_{L^p}^p= \int_\Rn \left|\frac{\chi_{\Omega_t}}{\Ln(\Omega_t)}-\frac{\chi_{\Omega_s}}{\Ln(\Omega_s)}\right|^p\,d\Ln 
	\geq \frac{\Ln(\Omega_t\triangle \Omega_s)}{M^p}.
	\end{align*}
	Therefore, if $|\dot {f}|^{L^p}(\cdot)\in L^1([0,T])$ is the $L^p$-metric derivative of $f_{(\cdot)}$, we have
	\begin{align}\label{rmk_ACLp_1}
	\Ln(\Omega_t\triangle \Omega_s)\leq M^p\left(\int_s^t|\dot{f}|^{L^p}(r)\,dr\right)^p.
	\end{align}
	Fix $\epsilon\ge 0$ arbitrarily. By standard properties of absolutely continuous functions, there exists a partition $t_0 = s \le t_1 \cdots \le t_N =t$ such that 
	\begin{align*}
	\int_{t_{j-1}}^{t_j}|\dot{f}|(r)\,dr \le \epsilon \quad \forall 1\leq j\leq N.
	\end{align*}
	Therefore, recalling the triangle inequality for the measure of the symmetric difference and using \eqref{rmk_ACLp_1} with $t_{j-1}$ and $t_j$, we obtain
	\begin{align*}
	\Ln(\Omega_t\triangle \Omega_s)&\leq \sum_{j=1}^N \Ln(\Omega_{t_j}\triangle \Omega_{t_{j-1}})\\
	& \leq M^p \epsilon^{p-1} \sum_{j=1}^N \int_{t_{j-1}}^{t_j}|\dot{f}|^{L^p}(r)\,dr\\
	& = M^p\epsilon^{p-1}\int_s^t|\dot{f}|^{L^p}(r)\,dr.
	\end{align*}
	Letting $\epsilon \downarrow 0$, we obtain exactly \eqref{rmk_ACLp}.
	
\end{remark}

	\subsection{Absolutely continuous curves in $\pl_q^p(\Rn)$ for $1<p\leq \infty$.}
\renewcommand{\s}{{\pl_q^p(\Rn)}}
\renewcommand{\d}{{\mathfrak{d}_q^p}}

Let us introduce the following notation. We say that a Borel function $F:[0,T]\times \Rn\to\R$ is \emph{$L^1$ in time and $L^p$ in space}, and write $F\in L^1_tL^p_x([0,T]\times \Rn)$ if for almost every $0\leq t\leq T$ the $t$-section $F_t:=F(t,\cdot)$ belongs to $L^p(\Rn)$ and
\begin{align*}
\int_0^T\norm{F_t}_{L^p}\,dt<\infty.
\end{align*}

\begin{theorem}\label{th_ACWqp_pNotOne}
	Let $1<p,q\leq \infty$. A curve $\mu_{(\cdot)}=f_{(\cdot)}\Ln:[0,T]\to \s$ is $\d$-absolutely continuous if and only if $\mu_{(\cdot)}$ is $W_q$-absolutely continuous and $f_{(\cdot)}$ is $L^p$-absolutely continuous. In particular, if $\mu_{(\cdot)}$ is $\d$-absolutely continuous, then there exists a Borel time-dependent vector field $(t,x)\mapsto v_t(x)$ that satisfies the following properties:
	\begin{enumerate}
	\item \label{th_ACPinftyp_0} $\int_0^T\norm{v_t}_{L^1(\mu_t)}\,dt<\infty$
	\item \label{th_ACPinftyp_i} the vector-valued function $(t,x)\mapsto v_t f_t(x)$ admits a weak $x$-divergence $\div(v f)\in L^1_tL^p_x([0,T]\times \Rn)$;
	\item \label{th_ACPinftyp_ii} the real-valued function $(t,x)\mapsto f_t(x)$ admits a weak $t$-derivative $\partial_t f\in L^1_tL^p_x([0,T]\times \Rn)$;
	\item \label{th_ACPinftyp_iibis} $\partial_t f + \div(vf) = 0$ in $L^1_tL^p_x([0,T]\times \Rn)$;
	\item \label{th_ACPinftyp_iiii} the function $t\mapsto \norm{v_t}_{L^q(\mu_t)}$ belongs to $L^1([0,T])$;
	\item \label{th_ACPinftyp_iv} $|\dot\mu|(t) = \norm{v_t}_{L^q(\mu_t)} + \norm{\div(v_t f_t)}_{L^p}$ for almost every $0\leq t\leq T$. 
	\end{enumerate}
	Conversely, if $\mu_{(\cdot)}=f_{(\cdot)}\Ln:[0,T]\to \s$ is narrowly continuous and there exists a Borel time-dependent vector field $(t,x)\mapsto v_t(x)$ such that properties \ref{th_ACPinftyp_0}--\ref{th_ACPinftyp_iiii} hold, then $\mu_{(\cdot)}$ is $\d$-absolutely continuous and 
	\begin{align}\label{E_0}
	|\dot \mu|(t)\leq \norm{v_t}_{L^q(\mu_t)} + \norm{\div(v_tf_t)}_{L^p}\quad\text{for almost every }0\leq t\leq T.
	\end{align} 
\end{theorem}

\begin{proof}
	If $\mu_{(\cdot)}=f_{(\cdot)}\Ln:[0,T]\to \s$ is $\d$-absolutely continuous, then 
	\begin{align*}
	{W_q(\mu_{t},\mu_s)}+ {\norm{f_{t}-f_s}}_{L^p} = \d(\mu_{t},\mu_s) \leq \int_s^t |\dot\mu|(r)\,dr\quad \forall 0\leq s \le t \leq T. 
	\end{align*}
	Therefore $\mu_{(\cdot)}$ is $W_q$-absolutely continuous and $f_{(\cdot)}$ is $L^p$-absolutely continuous. Viceversa, if  $\mu_{(\cdot)}$ is $W_q$-absolutely continuous and $f_{(\cdot)}$ is $L^p$-absolutely continuous, then 
	\begin{align*}
	\d(\mu_{t},\mu_s) = {W_q(\mu_{t},\mu_s)}+ {\norm{f_{t}-f_s}}_{L^p} \leq \int_s^t (|\dot\mu|^{W_q}(r) + |\dot{f}|^{L^p}(r))\,dr.
	\end{align*}
	This proves the first claim. Moreover, by the very definition of the metric derivative, it follows that if $\mu_{(\cdot)}$ is $\d$-absolutely continuous, then 
	\begin{align}\label{eq_th_metricDerivatives}
	|\dot \mu|(t) = |\dot \mu|^{W_q}(t) + |\dot{f}|^{L^p}(t)\quad \text{for almost every }0\leq t\leq T.
	\end{align}
	
	Let us now fix a $\d$-absolutely continuous curve $\mu_{(\cdot)}$. Recalling Lemma \ref{lem_timeInvariance}, it is enough to prove the statement for a curve $\mu_{(\cdot)}$ parametrized by arc-length, i.e. with $|\dot\mu|(t)=1$ at almost every $0\leq t\leq T$.
	 By virtue of Theorem \ref{th_ACinPq}, if $1<q<\infty$, or Corollary \ref{cor_ACinPinfty}, if $q=\infty$, there exists a Borel vector field $(t,x)\mapsto v_{t}(x)$ such that 
	\begin{gather}
	\int_0^T\int_\Rn (\partial_t\phi_t + \langle \nabla \phi_t,v_t \rangle) f_t\,d\Ln
	\,dt = 0 \quad \forall \phi\in \C^\infty_c((0,T)\times \Rn),\label{eq_thACPinftyp_1}\\
	\norm{v_t}_{L^q(\mu_t)} = |\dot\mu|^{W_q}(t) \leq  1\quad \text{for almost every }0\leq t \leq T.\label{eq_thACPinftyp_2}
	\end{gather}
	
	Taking into account Lemma \ref{lem_Bochner}, as $f_{(\cdot)}$ is $L^p$-absolutely continuous with essentially bounded $L^p$-metric derivative, there exists a function $\partialf_{(\cdot)}\in L^\infty([0,T];L^p(\Rn))$, if $1<p<\infty$, or $\partialf_{(\cdot)}\in L_{w^*}^\infty([0,T];L^\infty(\Rn))$, if $p=\infty$, such that
\begin{gather}
	\lim_{h\to 0}\int_{\Rn} \psi \frac{f_{t+h}-f_{t}}{h}\,d\Ln = \int_{\Rn} \psi \partialf_t\,d\Ln \,\,\forall \psi \in \C^0_0(\Rn)\,\text{a.e. }0\leq t\leq T,\label{eq_Dtft_1}\\
	 f_b-f_a = \int_{a}^b \partialf_t\,dt  \quad\text{in }L^p(\Rn),\,\forall 0\le a\leq b \le T,\label{eq_Dtft_2}\\
	\norm{\partialf_t}_{L^p}= |\dot {f}|(t)\leq 1 \quad \text{for a.e. }0\leq t\leq T.\label{eq_Dtft_3}
\end{gather}	 

	Fix a smooth compactly supported function $\phi\in \C^\infty_c((0,T)\times\Rn)$ and, for every $0\leq t \leq T$, denote by $\phi_t\in \C^\infty_c(\Rn)$ the function $x\mapsto \phi(t,x)$. Then
\begin{align*}
	\begin{split}
	\int_0^T\int_\Rn \partial_t \phi_t f_t\,d\Ln\,dt 
	&= \lim_{h\downarrow 0}\int_\Rn\int_0^T \frac{\phi_t - \phi_{t-h}}{h}f_t\,dt\,d\Ln\\
	&= \lim_{h\downarrow 0}\int_{\Rn}\frac{1}{h}\left(\int_0^T \phi_t f_t\,dt - \int_0^T \phi_{t-h}f_t\,dt\right)\,d\Ln.
	\end{split}
\end{align*}
	Since there exists $\delta \ge 0$ and a compact $K\subseteq \Rn$ such that $\spt \phi\subseteq [\delta, T-\delta]\times K$, we can apply the change of variables $t\mapsto t-h$ to the second integral to obtain
 \begin{align}\label{eq_ACPinftp_1}
	\int_0^T \int_\Rn \partial_t\phi_t f_t\,d\Ln\,dt = -\lim_{h\downarrow 0}\int_0^T\int_\Rn \phi_t\frac{f_{t+h}-f_t}{h}\,d\Ln\,dt.
\end{align}
From \eqref{eq_Dtft_1}, we deduce that
\begin{align*}
	\lim_{h\downarrow 0}\int_\Rn \phi_t \frac{\,f_{t+h}-f_t}{h}\,d\Ln = \int_\Rn \phi_t \partialf_t\,d\Ln \quad \text{for a.e. }0\leq t \leq T.
\end{align*}
Since
\begin{align*}
	\left|\int_\Rn \phi_t \frac{\,f_{t+h}-f_t}{h}\,d\Ln\right| \leq \frac{\norm{\phi}_{\C^0}}{h}\int_t^{t+h}\int_K |\partialf_s|\,d\Ln\,ds \leq \norm{\phi}_{\C^0}\left(\Ln(K)\right)^\frac{1}{p'},
\end{align*}
where $1\leq p'<\infty$ is the H\"older conjugate of $p$, by dominated convergence we deduce that
\begin{align}\label{eq_ACPinftp_2}
\lim_{h\downarrow 0} \int_0^T\int_\Rn\phi_t \frac{\,f_{t+h}-f_t}{h}\,d\Ln = \int_0^T \int_\Rn \phi_t \partialf_t\,d\Ln\,dt.
\end{align}
Combining \eqref{eq_thACPinftyp_1}, \eqref{eq_ACPinftp_1} and \eqref{eq_ACPinftp_2}, we obtain 
\begin{align}\label{eq_ACPinftyp_3}
	-\int_0^T\int_\Rn \langle \nabla \phi_t,v_t f_t\rangle\,d\Ln\,dt = \int_0^T \int_\Rn \partial_t\phi_t f_t\,d\Ln\,dt = - \int_0^T \int_\Rn \phi_t \partial f_t\,d\Ln\,dt.
\end{align}
Equation \eqref{eq_ACPinftyp_3}, together with \eqref{eq_Dtft_3}, proves \ref{th_ACPinftyp_i}, \ref{th_ACPinftyp_ii} and \ref{th_ACPinftyp_iibis}. Finally, \eqref{eq_th_metricDerivatives}, \eqref{eq_thACPinftyp_2}, and \eqref{eq_Dtft_3} imply \ref{th_ACPinftyp_iv}.

Suppose now $\mu_{(\cdot)}=f_{(\cdot)}\Ln:[0,T]\to \pl_q^p(\Rn)$ to be narrowly continuous and the existence of a Borel vector field $(t,x)\mapsto v_t(x)$ such that properties \ref{th_ACPinftyp_0}--\ref{th_ACPinftyp_iiii} hold. Then $v$ satisfies the integrability condition \eqref{def_contEqIntegrability} of Definition \ref{def_WeakSol} and, combining \ref{th_ACPinftyp_i}, \ref{th_ACPinftyp_ii} and \ref{th_ACPinftyp_iibis}, it immediately follows that $(\mu_{(\cdot)},v_{(\cdot)})$ is a weak solution of the continuity equation in $[0,T]$. Thus, from Theorem \ref{th_ACinPq} (resp. Theorem \ref{th_ineqA}, if $q=\infty$), we deduce that $\mu_{(\cdot)}$ is $W_q$-absolutely continuous and that 
\begin{align}\label{E_1}
|\dot\mu|^{W_q}(t)\leq \norm{v_t}_{L^q(\mu_t)} \quad \text{for almost every }0\leq t\leq T.
\end{align}

Fix $0\leq t\leq t+h\leq T$. We show that
\begin{align}\label{eqq_1}
f_{t+h} = f_t - \int_t^{t+h}\div(v_\tau f_\tau)\,d\tau\quad \text{in }L^p(\Rn),
\end{align}
where the integral in the right-hand side is the Bochner integral (resp. the weak$^*$ integral, if $p=\infty$) of $\tau \mapsto \div(v_\tau f_\tau)$. First, observe that -- thanks to \ref{th_ACPinftyp_i} and \ref{th_ACPinftyp_iiii} --  $\tau \mapsto \div(v_\tau f_\tau)$ is a Bochner measurable (resp. weakly$^*$ measurable, if $p=\infty$) curve $[0,T]\to L^p(\Rn)$, and it belongs to $L^1([0,T];L^p(\Rn))$ (resp. $L^1_{w^*}([0,T];L^\infty(\Rn))$, if $p=\infty$). Therefore, the right-hand side of \eqref{eqq_1} is well-defined.

We shall now prove the claim \eqref{eqq_1}. It will be enough to show that
\begin{align}\label{eqq_3}
\int_\Rn \left(f_{t+h} - f_t + \int_t^{t+h}\div(v_\tau f_\tau)\,d\tau\right)g\,d\Ln = 0\quad \forall g\in \C^\infty_c(\Rn).
\end{align}
Fix a smooth test function $g\in \C^\infty_c(\Rn)$. Then, recalling Remark \ref{rmk_NarrowContWeakSol}, the function $\tau\mapsto \int_\Rn f_\tau g\,d\Ln$ is absolutely continuous, and in particular
\begin{align*}
\int_\Rn f_{t+h}g\,d\Ln - \int_\Rn f_tg\,d\Ln = \int_t^{t+h}\int_\Rn \langle \nabla g,v_\tau\rangle\,d\Ln\,d\tau.
\end{align*} 
Therefore,
\begin{align}\label{eqq_2}
\begin{split}
\int_\Rn &\left(f_{t+h} - f_t + \int_t^{t+h}\div(v_\tau f_\tau)\,d\tau\right)g\,d\Ln 
= \\
&\quad\quad\int_0^{T}\int_\Rn \left(\langle \nabla g\,\chi_{[t,t+h]}(\tau),v_\tau\rangle + \div(v_\tau f_\tau)g\,\chi_{[t,t+h]}(\tau)\right)\,d\tau\,d\Ln.
\end{split}
\end{align}
Approximating $\chi_{[t,t+h]}$ with a smooth function that is compactly supported in $[0,T]$ and integrating by parts, we prove that the right-hand side of \eqref{eqq_2} is zero. By arbitrariness of the choice of $g$, \eqref{eqq_3} follows.

Using \eqref{eqq_1} and the triangle inequality for integrals, we can write
\begin{align*}
\norm{f_{t+h}-f_t}_{L^p}\leq \int_t^{t+h}\norm{\div(v_\tau f_\tau)}_{L^p}\,d\tau.
\end{align*}
Therefore, by arbitrariness of $t$ and $h$, and taking into account \ref{th_ACPinftyp_iiii}, $f_{(\cdot)}$ is $L^p$-absolutely continuous and 
\begin{align}\label{E_2}
|\dot f|^{L^p}(t)\leq \norm{\div(v_t f_t)}_{L^p}\quad \text{for almost every }0\leq t\leq T.
\end{align}

Combining \eqref{E_1} and \eqref{E_2}, we deduce that $\mu_{(\cdot)}$ is $\d$-absolutely continuous and that \eqref{E_0} holds true.

\end{proof}

\subsection{Absolutely continuous curves in $\pl_q^1(\Rn)$.}
\renewcommand{\s}{{\pl_q^1(\Rn)}}
\renewcommand{\d}{{\mathfrak{d}_q^1}}

We say that a Borel measure $M = \L^1\llcorner[0,T]\otimes (M_t)_{t}\in \mathcal{M}([0,T]\times \Rn)$ is \emph{$\mathcal{M}^1$ in time and of finite variation in space}, and write $M\in \mathcal{M}^1_t\text{\normalfont FV}_x([0,T]\times \Rn)$ if
\begin{align*}
\int_0^T\norm{M_t}_{\tv}\,dt<\infty.
\end{align*}

\begin{theorem}\label{th_AC_PLOne}
Let $1<q\leq \infty$.
	A curve $\mu_{(\cdot)}=f_{(\cdot)}\Ln:[0,T]\to \s$ is $\d$-absolutely continuous if and only if $\mu_{(\cdot)}$ is $W_q$-absolutely continuous and $f_{(\cdot)}$ is $L^1$-absolutely continuous. In particular, if $\mu_{(\cdot)}$ is $\d$-absolutely continuous, then there exists a Borel time-dependent vector field $(t,x)\mapsto v_t(x)$ that satisfies the following properties:
	\begin{enumerate}
	\item \label{th_ACPinfty1_0} $\int_0^T\norm{v_t}_{L^1(\mu_t)}\,dt<\infty$
	\item \label{th_ACPinfty1_i} the vector-valued function $(t,x)\mapsto v_t f_t(x)$ admits a weak $x$-divergence $\L^1\llcorner [0,T] \otimes(\div(v_t f_t))_t\in \mathcal{M}^1_t\text{\normalfont FV}_x([0,T]\times \Rn)$;
	\item \label{th_ACPinfty1_ii} the real-valued function $(t,x)\mapsto f_t(x)$ admits a weak $t$-derivative $\L^1\llcorner[0,T]\otimes(\partial_t f_t)_t\in \mathcal{M}^1_t\text{\normalfont FV}_x([0,T]\times \Rn)$;
	\item \label{th_ACPinfty1_iibis} $\L^1\llcorner[0,T]\otimes(\partial_t f_t+\div(v_t f_t))_t = 0$ in $\mathcal{M}^1_t\text{\normalfont FV}_x([0,T]\times \Rn)$;
	\item \label{th_ACPinfty1_iiii} the function $t\mapsto \norm{v_t}_{L^q(\mu_t)}$ belongs to $L^1([0,T])$;
	\item \label{th_ACPinfty1_iv} $|\dot\mu|(t) = \norm{v_t}_{L^q(\mu_t)} + \norm{\div(v_t f_t)}_{\tv}$ for almost every $0\leq t\leq T$. 
	\end{enumerate}
	Conversely, if $\mu_{(\cdot)}=f_{(\cdot)}\Ln:[0,T]\to \s$ is narrowly continuous and there exists a Borel time-dependent vector field $(t,x)\mapsto v_t(x)$ such that properties \ref{th_ACPinfty1_0}--\ref{th_ACPinfty1_iiii} hold, then $\mu_{(\cdot)}$ is $\d$-absolutely continuous and 
	\begin{align}\label{F_0}
	|\dot \mu|(t)\leq \norm{v_t}_{L^q(\mu_t)} + \norm{\div(v_tf_t)}_{\tv}\quad\text{for almost every }0\leq t\leq T.
	\end{align} 
\end{theorem}

\begin{proof}
	The first claim is proved exactly as in Theorem \ref{th_ACWqp_pNotOne}. 
	
	Fix a curve $\mu_{(\cdot)}=f_{(\cdot)}\Ln\in \ac([0,T];\pl_q^1(\Rn))$.
	Consider the isometric embedding of $(L^1(\Rn),\norm{\cdot}_{L^1})$ into the space of signed real measures $(\mathcal{M}(\Rn),\norm{\cdot}_{\tv})\simeq (\C^0_b(\Rn),\norm{\cdot}_{\C^0})^*$ given by $f\hookrightarrow f\Ln$. 
	Therefore, the curve $f_{(\cdot)}\in \ac^\infty([0,T];L^1(\Rn))$ is identified (via the aforementioned embedding) with the curve $f_{(\cdot)}\Ln\in \ac^\infty([0,T];\mathcal{M}(\Rn))$. Hence, by virtue of Lemma \ref{lem_Bochner}, there exists a curve $\partialf_{(\cdot)}\in L^\infty_{w^*}([0,T];\mathcal{M}(\Rn))$ such that
	\begin{gather*}
	\lim_{h\to 0}\int_{\Rn} \psi \frac{f_{t+h}-f_{t}}{h}\,d\Ln = \int_{\Rn} \psi\, d\partialf_t \quad\forall \psi \in \C^0_0(\Rn)\,\text{ for a.e. }0\leq t\leq T,\\
	 \int_\Rn \psi\,(f_b-f_a)\,d\Ln = \int_{a}^b\int_\Rn \psi\, d\partialf_t\,dt  \quad\forall \psi \in \C^0_0(\Rn)\,\forall 0\le a\leq b \le T,\\
	\norm{\partialf_t}_{\tv}= |\dot {f}|(t)\leq 1 \quad \text{for a.e. }0\leq t\leq T.
\end{gather*}
	
	Fix a test function $\phi\in \C^\infty_c([0,T]\times \Rn)$. Arguing as in the proof of Theorem \ref{th_ACWqp_pNotOne}, we obtain the equality
	\begin{align*}
	\int_0^T\int_\Rn \partial_t\phi_t f_t\,d\Ln\,dt = -\lim_{h\downarrow 0}\int_0^T\int_\Rn\phi_t\frac{f_{t+h}-f_t}{h}\,d\Ln\,dt,
\end{align*}	 	 
and the limit
\begin{align*}
	\lim_{h\downarrow 0} \int_\Rn\phi_t\frac{f_{t+h}-f_t}{h}\,d\Ln = \int_\Rn\phi_t\,d\partialf_t\quad \text{for almost every }0\leq t\leq T.
\end{align*}
To justify the exchange of limit and integral, it is enough to observe that
\begin{align*}
\left|\int_\Rn \phi_t \frac{\,f_{t+h}-f_t}{h}\,d\Ln\right| \leq \frac{\norm{\phi}_{\C^0}}{h}\int_t^{t+h}\norm{\partialf_s}_{\tv}\,ds \leq \norm{\phi}_{\C^0}
\end{align*}
and to use the dominated convergence theorem. Thus, if $\mu_{(\cdot)}$ is $\d$-absolutely continuous and $v_{(\cdot)}$ is the velocity field of $\mu_{(\cdot)}$ given by Theorem \ref{th_ACinPq}, for $1<q<\infty$, or by Corollary \ref{cor_ACinPinfty}, if $q=\infty$, then
\begin{align*}
-\int_0^T\int_\Rn \langle \nabla \phi_t,v_t\rangle f_t\,d\Ln\,dt = \int_0^T\int_\Rn \partial_t\phi_t\,f_t\,d\Ln = -\int_0^T\int_\Rn \phi_t\,d\partialf_t\,dt,
\end{align*}
and properties \ref{th_ACPinfty1_0}--\ref{th_ACPinfty1_iv} are proved.

Suppose now $\mu_{(\cdot)}=f_{(\cdot)}\Ln:[0,T]\to \s$ is narrowly continuous and there exists a Borel time-dependent vector field $(t,x)\mapsto v_t(x)$ such that properties  \ref{th_ACPinfty1_0}--\ref{th_ACPinfty1_iiii} hold. Then, arguing as in Theorem \ref{th_ACWqp_pNotOne} one proves that $(\mu_{(\cdot)},v_{(\cdot)})$ is a weak solution of the continuity equation in $[0,T]$. Therefore $\mu_{(\cdot)}$ is $W_q$-absolutely continuous and the estimate 
\begin{align*}
|\dot \mu|^{W_q}(t)\leq \norm{v}_{L^q(\mu_t)}
\end{align*}
holds for almost every $0\leq t\leq T$. On the other hand, using the absolute continuity of the map $\tau\mapsto \int_\Rn g\,d\mu_\tau$ given by Remark \ref{rmk_NarrowContWeakSol}, for every $g\in \C^\infty_c(\Rn)$, together with a smooth compactly supported approximation of $\chi_{[t,t+h]}$, one proves that
\begin{align*}
\mu_{t+h} = \mu_t + \int_t^{t+h}\div(v_\tau f_\tau)\,d\tau\quad \text{in }\mathcal{M}(\Rn)
\end{align*}
for all $0\leq t\leq t+h\leq T$. Therefore, 
\begin{align*}
\norm{f_{s}-f_t}_{L^1} = \norm{\mu_{s}-\mu_t}_{\tv}\leq \int_t^{s}\norm{\div(v_\tau f_\tau)}_{\tv}\,d\tau\quad \forall 0\leq t\leq s\leq T.
\end{align*}
Hence $f_{(\cdot)}$ is $L^1$-absolutely continuous and \eqref{F_0} holds.

\end{proof}

\renewcommand{\d}{{\mathfrak{d}_q^p}}
\subsection{Examples of $\d$-AC curves.}

Theorem \ref{th_ACWqp_pNotOne} and Theorem \ref{th_AC_PLOne} give a complete characterization of $\d$-absolutely continuous curves in terms of solution of the continuity equations for velocity fields $v_{(\cdot)}$ that satisfy a Sobolev-like condition involving the $L^q(\mu_t)$-norm of $v_t$ and the $L^p(\Rn)$-norm of $\div(vf)$. However, Theorems  \ref{th_ACWqp_pNotOne} and \ref{th_AC_PLOne} do not guarantee actual existence of such fields, nor of $\mathfrak{d}_q^p$-absolutely continuous curves.

 We shall now give some examples of $\d$-absolutely continuous curves.

\renewcommand{\d}{\mathfrak{d}_\infty^\infty}
\begin{example}[interpolation of densities bounded from below]
	Let $\mu =f\Ln,\nu = g\Ln\in \pl^p_\infty(\Rn)$, for some $n < p\leq \infty$. Suppose that $K:=\spt\,f =\spt \,g$ is a smooth and connected domain and the existence of a constant $0<c$ such that $\min\{f,g\}\geq c$ almost everywhere in $K$. Then there exists a $\mathfrak{d}^p_\infty$-absolutely continuous curve $\mu_{(\cdot)}$ that connects $\mu$ to $\nu$. 
	
	To prove this claim, define the interpolation function $f_t:\Rn\to [0,\infty)$
	\begin{align*}
	f_t(x):=(1-t)f(x) + tg(x)\quad \forall 0\leq t\leq 1.
	\end{align*}
	Observe that $\mu_t:=f_t\Ln$ is a probability measure and 
	\begin{align*}
	\norm{f_{s}-f_t}_{L^p}\leq (s-t) \norm{f-g}_{L^p} \quad \forall 0\leq t\leq s\leq 1.
	\end{align*}
	Therefore 
	the curve of densities $f_{(\cdot)}$ is $L^p$-absolutely continuous.
	The function $t\mapsto f_t(x)$ is differentiable for every $x$ and 
	\begin{align*}
	\partial_t f_t(x) = g(x)-f(x)\quad \forall x\in \Rn,\,0\leq t\leq 1.
	\end{align*}
	Therefore, recalling that $\mu$ and $\nu$ are probability measures,
	\begin{align*}
	\int_\Rn \partial_tf_t\,d\Ln = 0.
	\end{align*}
	By Bogovski\u{\i} theorem (see e.g. Lemma III.3.1 in \cite{Galdi2011}), if $n<p<\infty$ there exists $w\in W^{1,p}_0(K;\Rn)$ such that 
\begin{gather}\label{eq_Bog}
\begin{gathered}
\div \,w = -\partial_t f_t\quad \text{in }K,\\
\norm{w}_{W^{1,p}(K)}< C\norm{\partial_tf_t}_{L^p(K)} = C\norm{f-g}_{L^p},
\end{gathered}
\end{gather}
for some constant $C$ that depends solely on $K$ and $p$. If $p=\infty$, then fix $n<\overline p<\infty$ and use Bogovski\u{\i} theorem to find $w\in W^{1,\overline{p}}_0(K;\Rn)$ such that properties \eqref{eq_Bog} hold with $\overline{p}$ in place of $p$.  In any case Sobolev embedding gives $w\in L^\infty(K;\Rn)$. 
Define  $v_t=w/f_t$ on $K$. Since both $f$ and $g$ are essentially uniformly bounded from below in $K$ by a positive constant, then so it is $f_t$. This implies the field $v_t:=w/f_t$ is well-defined almost everywhere in $K$ and
\begin{align*}
	\norm{v_t}_{L^\infty(K)}\leq \frac{\norm{w_t}_{L^\infty}}{c}\leq \hat{C}\norm{f-g}_{L^p}\quad \forall 0\leq t\leq 1,
	\end{align*}
	where $\hat{C}:= C/c$.
	
	On the other hand, integrating by parts, we deduce that
	\begin{align*}
	\int_0^1\int_\Rn \langle \nabla \psi(t,x),w_t(x)\rangle dx\,dt 
	&= \int_0^1\int_\Rn \psi(t,x)\partial_tf_t(x)\,dx\,dt\\
	&= -\int_0^1\int_\Rn \partial_t\psi(t,x)f_t(x)\,dx\,dt
	\end{align*}
	for all $\psi\in \C^\infty_c((0,1)\times \Rn)$. 
Since $(\mu_{(\cdot)},v_{(\cdot)})$ is a narrowly continuous solution of the continuity equation, Theorem \ref{th_ineqA} ensures the $W_\infty$-absolute continuity of $\mu_{(\cdot)}$. Therefore $\mu_{(\cdot)}$ is $\mathfrak{d}^p_\infty$-absolutely continuous (and therefore $\mu_{(\cdot)}$ is also $\mathfrak{d}^p_q$-absolutely continuous for every $1 < q\leq \infty$).
	
\end{example}

\begin{example}[Translations and dilations of compact Lipschitz densities]\label{examp_TranDil}
	Let $\mu=f\Ln\in \pl^\infty_\infty(\Rn)$ be such that $f\in \Lip_c(\Rn)$. Let $v\in \Rn$ be a vector and $M>0$. The translation curve of $\mu$ by $v$ and the dilation curve by $M$ of $\mu$, that are respectively defined as
	\begin{align*}
	\hat{\mu}_t^v:=f(\cdot -tv)\Ln \quad\text{and}\quad \tilde{\mu}^M_t:=\frac{f(\cdot/(1-t+tM))}{(1-t+t M)^n}\Ln\quad \forall 0\leq t\leq 1,
\end{align*}	  
are $\d$-absolutely continuous curves in $\pl^\infty_\infty(\Rn)$. 
Indeed, a velocity-field for $\hat\mu_{(\cdot)}$ is given by the constant field $\equiv v$, therefore $\norm{v_t}_{L^\infty(\mu)}\equiv |v|$ and the map $t\mapsto f(x-tv)$ is Lipschitz, with Lipschitz constant not greater than $\Lip(f)|v|$. Therefore, by Theorem \ref{th_ACWqp_pNotOne} it follows that the translation $\hat\mu_{(\cdot)}^v$ is $\d$-absolutely continuous and
\begin{align*}
|\dot {\hat{\mu}}^v|(t)\leq (\Lip(f) + 1) |v|\quad \text{for almost every }0\leq t\leq 1.
\end{align*}

Similarly, it is easy to verify that a velocity field for $\tilde\mu_{(\cdot)}$ is given by
\begin{align*}
v_t(x):= \frac{(M-1)}{1+t(M-1)} x.
\end{align*}
Clearly $\norm{v_t}_{L^\infty(\mu_t)}\leq R|M-1|<\infty$, if $(1+M)\spt\,\mu\subseteq B_R$, and the function $t\mapsto f(x/(1-t+tM))(1-t+tM)^{-n}$ is Lipschitz, with Lipschitz constant uniformly bounded in $x$. This proves the $\d$-absolute continuity of $\tilde\mu_{(\cdot)}$.

\end{example}

\section[final remarks and open questions]{Final remarks and open questions}\label{sec_Rmk}

We conclude the paper by collecting some natural questions and possible directions for further research. We hope that these problems may be of interest to other researchers and may contribute to the further development of the metric framework introduced in this paper.

This paper shows that the minimizing movement scheme for the isoperimetric functional admits generalized minimizing movements which are absolutely continuous in $\pl_q^p(\Rn)$. We also obtained a characterization of absolutely continuous curves in $\pl_q^p(\Rn)$ in terms of the continuity equation and suitable Eulerian estimates. 

These results do not imply, however, that $\pl_q^p(\Rn)$ contains many non-constant absolutely continuous curves. The characterization theorem is a criterion, not an existence result for curves with prescribed endpoints. In principle, it could happen that the metric $\mathfrak{d}_q^p$ is so restrictive that the only absolutely continuous curves starting from certain measures are constant curves. In such a situation, the GMMs constructed by the scheme with these initial data would be dynamically trivial. It is therefore natural to ask whether $\pl_q^p(\Rn)$ is connected by absolutely continuous arcs in the following sense.

\begin{problem}
Let $1<q\leq\infty$ and $1\leq p\leq\infty$. Given two arbitrary $\mu,\nu\in \pl_q^p(\mathbb R^n)$, does a $\mathfrak{d}_q^p$-absolutely continuous curve $\mu_{(\cdot)}:[0,1]\to \pl_q^p(\Rn)$ such that $\mu_0=\mu$ and $\mu_1=\nu$ always exists?
\end{problem}

We finally mention that Problem 1 has been solved, in the case $q=\infty$, in a work in progress by the author \cite{AldrigoInPreparation}. The general case with $1<q<\infty$ remains open.

Another natural question concerns the stationary points of the minimizing movement scheme of $\Isop$. Let $\mu\in \pl_\infty^\infty(\mathbb R^n)$. We say
that $\mu$ is \emph{stationary for the isoperimetric minimizing movement} if every GMM for $\operatorname{Isop}$ starting from $\mu$ is constant. The collection of these probability measures is denoted by $\Crit(\Isop)$.

Such measures should be regarded as critical points of the isoperimetric functional, in the variational sense given by the minimizing movement scheme. They play the role of configurations for which the discrete implicit Euler scheme detects no descending direction.

A basic class of stationary configurations is given by normalized Lebesgue measures on Euclidean balls -- as these are in fact the only \emph{global minima of $\Isop$}. More generally, one is led to expect that a particular subfamily of weighted sums of uniform measures on pairwise disjoint closed balls are stationary for the isoperimetric minimizing movement.

\begin{problem}
Characterize all the measures $\mu\in \pl_\infty^\infty(\Rn)$ that are stationary for the isoperimetric minimizing movement. In particular, is it true that a measure $\mu\in\pl_\infty^\infty(\mathbb R^n)$ is stationary for the isoperimetric minimizing movement if and only if it is of the form
\begin{align*}
    \mu_{\mathbf{r}}:=\sum_{j=1}^N \frac{1}{R\omega_nr_j^{n-1}}\Ln\llcorner \overline{B}^j, \quad \mathbf{r}=(r_j)_{1\leq j\leq N}\subseteq (0,\infty),\quad R:=\sum_{j=1}^N r_j,
\end{align*}
for some integer $N\geq 1$, where $(\overline B^j)_{1\leq j\leq N}$ are 
pairwise disjoint closed Euclidean balls with radii $(r_j)_{1\leq j\leq N}$?
\end{problem}

Finally, we ask for which initial data the minimizing movements for $\Isop$ always converge, as $t\to\infty$, to a stationary configuration, and for which of these initial data the limiting configuration is uniquely determined by the initial one. 

\begin{problem}
For which $\mu\in\pl_\infty^\infty(\Rn)$ every curve $\mu_{(\cdot)}\in \GMM(\Isop,\mu)$ has finite $\mathfrak{d}_\infty^\infty$-length, i.e.
\begin{align*}
\int_0^\infty |\dot \mu|(t)\,dt<\infty,
\end{align*} 
and there exists $\mu_\infty\in \Crit(\Isop)$ such that $\mu_t\to \mu_\infty$ in $\pl_\infty^\infty(\Rn)$ as $t\to\infty$? For which $\mu\in \pl_\infty^\infty(\Rn)$ is $\mu_\infty$ unique?
\end{problem}

\begin{appendices}

\titleformat{\section}
    {\large \centering \bfseries} 
    {Appendix \thesection.} 
    {5pt} 
    {} 

\section[dynamic formulations of the $\infty$-wasserstein metric]{Dynamic formulations of the $\infty$-Wasserstein metric}
\renewcommand{\d}{\text{\normalfont
d}}
\newcommand{\A}{\mathscr{A}_\infty}
\newcommand{\CC}{{\C^0([0,1];\Sp)}}
\newcommand{\AC}{{\ac([0,1];\Sp)}}

\subsection{A Benamou--Brenier formula.}

\begin{lemma}\label{lem_tech}
	Let $X,Y$ be two Polish space and let $\alpha\in \P(X\times Y)$. Denote by $\pi_X:X\times Y\to X$ be the canonical projection onto $X$ and let $\beta \otimes (\gamma_x)_x$ be a Borel disintegration of $\alpha$ with respect to $\pi_X$. Then, for every Borel-measurable function $f:X\times Y\to \R$:
	\begin{enumerate}
	\item \label{lem_tech_i} the function $y\mapsto f_x(y):=f(x,y)$ is Borel-measurable in $Y$ for all $x\in X$,
	\item \label{lem_tech_ii} the function $x\mapsto g_p(x):=\norm{f_x}_{L^p(\gamma_x)}$ is Borel-measurable in $X$ for all $1\leq p\leq \infty$, and
	\item \label{lem_tech_iii} $\norm{g_p}_{L^p(\beta)} = \norm{f}_{L^p(\alpha)}$ for all $1\leq p\leq \infty$.
	\end{enumerate}
\end{lemma}

\begin{theorem}\label{th_BenBre}
	For every couple of probability measures $\mu,\nu\in \P_\infty(\Rn)$, 
	\begin{align*}
	W_\infty(\mu,\nu) = \min\left\{\norm{v}_{L^\infty(\tilde{\mu})}: \begin{matrix}
	\tilde{\mu} :=\L^1\llcorner[0,1] \otimes(\mu_t)_t\hfill\\
	\partial_t\mu_t + \div(v_t\mu_t) = 0\text{ on }[0,1]\hfill\\
	\mu_{(\cdot)}\text{ narrowly continuous}\hfill\\
	\mu_0 = \mu\text{ and } \mu_1=\nu\hfill\\
	\end{matrix}
	\right\}.
	\end{align*}
\end{theorem}

\begin{proof}
Fix two arbitrary probability measures $\mu,\nu\in \P_\infty(\Rn)$.
	We begin by constructing a solution $(\mu_{(\cdot)},v_{(\cdot)})$ of the continuity equation in $[0,1]$ such that $\mu_0=\mu$, $\mu_1=\nu$ and with 
	\begin{align}\label{A_1}
	\norm{v}_{L^\infty(\tilde{\mu})}\leq W_\infty(\mu,\nu).
	\end{align}
	Let $\gamma\in \Gamma_\infty(\mu,\nu)$ be optimal and define the functions $e_t,V:\Rn\times\Rn\to \Rn$ by setting $V(x,y):=y-x$ and $e_t(x,y):=(1-t)x+ty$ for every $x,y\in \Rn$ and every $0\leq t\leq 1$. Consider now the measures $\mu_t:=(e_t)_\#\gamma\in \P(\Rn)$ and $E_t:=(e_t)_\#(V\gamma)\in \mathcal{M}(\Rn;\Rn)$ for every $0\leq t\leq 1$.
	Arguing as in the proof of Theorem 17.2 of \cite{ambrosio2021} and using Riesz representation theorem for the dual of $L^1(\mu_t)$, together with some standard measurability technicalities, one construct a velocity field $v_{(\cdot)}$ for $\mu_{(\cdot)}$ that satisfies
	\begin{align*}
	\norm{v_t}_{L^\infty(\mu_t)}\leq W_\infty(\mu,\nu)
	\end{align*}
	for almost every $0\leq t\leq T$. Therefore, \eqref{A_1} follows from Lemma \ref{lem_tech}.

	To prove the other inequality, we show that for every narrowly continuous solution $(\mu_{(\cdot)},v_{(\cdot)})$ of the continuity equation on $[0,1]$ such that $\mu_0=\mu$ and $\mu_1=\nu$ we have
\begin{align}\label{th_eq_BenBre_1}
	W_\infty(\mu,\nu)\leq \norm{v}_{L^\infty(\tilde{\mu})}.
\end{align}	
	 Let $\{\rho_\epsilon:\epsilon >0\}$ be a family of strictly positive and rapidly decreasing mollifiers, and let $(\mu^\epsilon_{(\cdot)},v^\epsilon_{(\cdot)})$ be the approximation of $(\mu_{(\cdot)},v_{(\cdot)})$ defined in \eqref{lem_def_approx} of Lemma \ref{Lem_Approx}. By Proposition 8.1.8 of \cite{ambrosio2005gradient}, it follows that $\mu_t^\epsilon = (T^\epsilon_t)_\#\mu_0^\epsilon$ for every $0\leq t\leq 1$, where $(t,x)\mapsto T_t^\epsilon(x)$ is the flow-map of the time-dependent vector field $v_{(\cdot)}^\epsilon$. Therefore, setting $\gamma^\epsilon:=(\id\times T_1^\epsilon)_\#\mu_0^\epsilon\in \Gamma(\mu_0^\epsilon,\mu_1^\epsilon)$, we have
	\begin{align*}
	|x-y| = |x-T_1^\epsilon(x)| \leq \int_0^1 |v_t^\epsilon(T_t^\epsilon(x))|\,dt\quad \gamma-\text{a.e. }(x,y)\in \Rn\times \Rn.
	\end{align*}
	Therefore, recalling \eqref{lem_approx_2} of Lemma \ref{Lem_Approx}, we infer
	\begin{align*}
	W_\infty(\mu_0^\epsilon,\mu_1^\epsilon) &\leq \norm{\int_0^1|v_t^\epsilon(T_t^\epsilon(x))|\,dt}_{L^\infty(\mu_0)} \leq \int_0^1 \norm{v_t^\epsilon}_{L^\infty(\mu_t^\epsilon)}\,dt \\
	&\leq \norm{t\mapsto \norm{v_t}_{L^\infty(\mu_t)}}_{L^\infty([0,1])}.
	\end{align*}
	Finally, combining the narrow convergence of $\mu_0^\epsilon$ and $\mu_1^\epsilon$ to $\mu_0=\mu$ and $\mu_1=\nu$ respectively guaranteed by \eqref{lem_approx_4} of Lemma \ref{Lem_Approx}, Lemma \ref{lem_WinftyIsLSC} and Lemma \ref{lem_tech}, \eqref{th_eq_BenBre_1} follows.

\end{proof}

\subsection{Probabilistic representation of $W_\infty$-AC curves.} 
\renewcommand{\CC}{{\C^0([0,T];\Rn)}}
\renewcommand{\AC}{{\ac([0,T];\Rn)}}

\begin{theorem}\label{th_ProbRep}
	Let $\mu_{(\cdot)}:[0,T]\to \P(\Rn)$ be a narrowly continuous curve and let $v_{(\cdot)}$ be a velocity-field for $\mu_{(\cdot)}$. Denote by $\tilde{\mu}$ the measure $\L^1\llcorner[0,T]\otimes(\mu_t)_t$ and suppose that $\norm{v}_{L^q(\tilde{\mu})}<\infty$ for some $1<q<\infty$. Then there exists a probability measure $\eta\in \P(\CC)$ such that
\begin{enumerate}
	\item \label{th_ProbRep_i} $\eta$ is concentrated on the family $\{\omega \in \AC: \dot\omega(t) = v_t(\omega(t))\,\ae\}$, and
	\item \label{th_ProbRep_ii} $(e_t)_\#\eta = \mu_t$ for every $0\leq t\leq T$.
\end{enumerate}	
\end{theorem}

We refer to Theorem 8.2.1 in \cite{ambrosio2005gradient} for a proof of the above statement. 

\begin{theorem}\label{cor_ProbabilisticRepresentation}
	A narrowly continuous curve $\mu_{(\cdot)}:[0,T]\to \P_\infty(\Rn)$ is absolutely continuous if and only if there exists a Borel time-dependent vector field $(t,x)\mapsto v_t(x)$ and $\eta\in \P(\CC)$ such that 
	\begin{enumerate}
	\item \label{cor_ProbRep_i} $\eta$ is concentrated on the family $\left\{\omega\in \AC : \dot \omega = v\circ\omega \text{ a.e. }\right\}$,
	\item \label{cor_ProbRep_ii} $\int_0^T \norm{v_t}_{L^\infty(\mu_t)}\,dt <\infty$, and
	\item \label{cor_ProbRep_iii} $(e_t)_\#\eta = \mu_t$ for every $0\leq t \leq T$.
	\end{enumerate}
\end{theorem}

\begin{proof}
	Suppose the existence of $v$ and $\eta$ that satisfy \ref{cor_ProbRep_i}, \ref{cor_ProbRep_ii} and \ref{cor_ProbRep_iii}. Arguing as in the first part of the proof of Theorem 8.2.1 of \cite{ambrosio2005gradient} and using Corollary \ref{cor_ACinPinfty}, we deduce that $\mu_{(\cdot)}\in \ac([0,T];\P_\infty(\Rn)))$. 
	
	Suppose now that  $\mu_{(\cdot)}\in \ac([0,T];\P_\infty(\Rn))$ is parametrized by arc-length. By virtue of Corollary \ref{cor_ACinPinfty}, $\mu_{(\cdot)}$ admits a velocity field $v_{(\cdot)}$ with $\norm{v_t}_{L^\infty(\mu_t)}=1$ for almost every $0\leq t\leq T$. Since $\tilde{\mu}:=\L^1\llcorner[0,T]\otimes(\mu_t)_t$ is a finite measure on $[0,T]\times \Rn$ and $v\in L^\infty(\tilde{\mu})$, then Theorem \ref{th_ProbRep} applies and \ref{cor_ProbRep_i}, \ref{cor_ProbRep_ii} and \ref{cor_ProbRep_iii} follow.
	
	If $\mu_{(\cdot)}\in \ac([0,T];\P_\infty(\Rn))$ is a general curve, then let $S:=\norm{\,|\dot\mu|^{W_\infty}}_{L^1([0,T])}$, $\sigma:[0,T]\to [0,S]$ the monotone reparametrization
	\begin{align*}
	\sigma(t):=\int_0^t|\dot \mu|(r)\,dr,
	\end{align*}
	 $\tau:[0,S]\to[0,T]$ its pseudo-inverse
	\begin{align*}
	\tau(s):=\min\{t\in [0,T]: \sigma(t) = s\},
\end{align*}		 
	 and $\nu_{(\cdot)}:= \mu_{\tau(\cdot)}:[0,S]\to \P_\infty(\Rn)$ be its arc-length  reparametrization. Thanks to the above arguments, there exist a Borel time-dependent vector field $w_{(\cdot)}$ and a probability measure in $\gamma\in \P(\C^0([0,S];\Rn))$  that satisfy the properties \ref{cor_ProbRep_i}, \ref{cor_ProbRep_ii} and \ref{cor_ProbRep_iii} for the curve $\nu_{(\cdot)}$. Let $\Phi_\sigma:\C^0([0,S];\Rn)\to \C^0([0,T];\Rn)$ be the function $\Phi_\sigma(\omega):=\omega\circ \sigma$ and define $\eta:=(\Phi_\sigma)_\#\gamma$. The function $\Phi_\sigma$ is injective. 
	 This implies that
	 \begin{align}\label{eq_cor_ProbRep_1}
	 \eta(\Phi_\sigma(\{\omega\in \ac([0,S];\Rn):\dot\omega(s)=w_s(\omega(s)) \text{ for a.e.  }0\leq s\leq S\})) = 1. 
	 \end{align}
	 Since $\sigma$ is monotone, then for every absolutely continuous curve $\omega\in \ac([0,S];\Rn)$ the curve $\Phi_\sigma(\omega)$ belongs to $\ac([0,T];\Rn)$ and its derivative is
	 \begin{align}\label{eq_cor_ProbRep_2}
	 \frac{d}{dt}\Phi_\sigma(\omega)(t) = \sigma'(t)\dot\omega(\sigma(t))\quad \text{for a.e. }0\leq t\leq T.
	 \end{align}
	 Therefore, defining $(t,x)\mapsto v_t(x):=\sigma'(t)w_\sigma(t)(x)$, and combining \eqref{eq_cor_ProbRep_1} and \eqref{eq_cor_ProbRep_2}, we deduce that $\eta$ is concentrated on the family
	 \begin{gather*}
	 \left\{\omega\in \ac([0,T];\Rn): \dot \omega(t) = v_t(\omega(t))\,\text{ for a.e. }0\leq t\leq T\right\}.
	 \end{gather*}
	 Moreover, as
	 \begin{align*}
	 W_\infty(\mu_t,\nu_{\sigma(t)}) = W_\infty(\mu_t,\mu_{\tau(\sigma(t))}) = \int_{\tau(\sigma(t))}^t |\dot \mu|(r)\,dr = \sigma(t)-\sigma(\tau(\sigma(t))) =0,
	 \end{align*}
	 holds for every $0\leq t\leq T$, 
	 we deduce that
	 \begin{align*}
	 \mu_t = \nu_{\sigma(t)} = (e_{\sigma(t)})_\#\gamma  = (e_t\circ \Phi_\sigma)_\#\gamma = (e_t)_\#\eta\quad \forall 0\leq t\leq T.
	 \end{align*}
	 Finally, using the substitution $s=\sigma(t)$,
	 \begin{align*}
	 \int_0^T\norm{v_t}_{L^\infty(\mu_t)}\,dt = \int_0^T\sigma'(t)\norm{w_\sigma(t)}_{L^\infty(\nu_{\sigma(t)})}\,dt = \int_0^S \norm{w_s}_{L^\infty(\nu_s)}\,ds <\infty.
	 \end{align*}
	 
\end{proof}

\subsection{Minimization of the $\infty$-Action.}
\renewcommand{\d}{{\normalfont\texttt{d}}}
\renewcommand{\CC}{\C^0([0,1];\mathscr{S})}

To simplify the notation, for the rest of the subsection we suppose $T=1$. Let $(\Sp,\d)$ be a general metric space.

\begin{definition}[Constant speed geodesics]
	A curve $\omega\in \ac([0,1];\Sp)$ is a \emph{constant speed geodesic}, and we write $\omega \in \csg(\Sp)$ if its metric derivative $|\dot\omega|$ is constantly equal to $\d(\omega(0),\omega(1))$ almost everywhere in $[0,1]$.  
\end{definition}

\begin{definition}[$\infty$-Action]
	The \emph{$\infty$-Action} is the function $\A:\CC\to [0,\infty]$ defined by
	\begin{align*}
	\A(\omega):=
	\begin{cases}
	\norm{\,|\dot\omega|\,}_{L^\infty([0,1])}&,\text{ if }\omega \in \ac([0,1];\Sp)\\
	\infty &,\text{ otherwise}
	\end{cases}.
	\end{align*}
\end{definition}

\begin{lemma}\label{lem_Ainfty}
	For every continuous curve $\omega\in \CC$ the inequality $\d(\omega(0),\omega(1))\leq \A(\omega)$ holds. Moreover, equality holds if and only if $\omega \in \csg(\Sp)$.
\end{lemma}

\begin{proof}
	The statement is trivial if $\omega \not \in \ac([0,1];\Sp)$. Suppose now $\omega\in \AC$. Then
	\begin{align*}
	\d(\omega(0),\omega(1))\leq \int_0^1|\dot \omega|(t)\,dt\leq \A(\omega),
	\end{align*}
	and equality holds if and only if $|\dot\omega|\equiv \A(\omega) = \d(\omega(0),\omega(1))$ almost everywhere in $[0,1]$.
	
\end{proof}

\begin{theorem}\label{th_DynamicFormulation}
	If $(\Sp,\d)$ is a Polish geodesic space and for every $0\leq t\leq 1$ denote by $e_t:\CC\to \Sp$ the $t$-evaluation $e_t(\omega):=\omega(t)$. Then, for every $\mu,\nu\in \P_\infty(\Sp)$ we have
	\begin{align}\label{th_DynForm_claim}
	W_\infty(\mu,\nu) = \min\left\{\norm{\A}_{L^\infty(\eta)}: \begin{matrix} \eta\in \P(\CC),\hfill \\ (e_0)_\#\eta = \mu,\hfill \\(e_1)_\#\eta =\nu\hfill \end{matrix}\right\}.
	\end{align}
	In particular, $\eta$ is a minimizer of if and only if 
	\begin{enumerate}
	\item \label{th_DynForm_i} $\norm{\A}_{L^\infty(\eta)} = \norm{\d(\cdot,\cdot)}_{L^\infty{(e_0\times e_1)_\#\eta}}$, and
	\item $(e_0\times e_1)_\#\eta \in \Gamma_\infty(\mu,\nu)$.
\end{enumerate}
	Moreover there always exists a minimizer of $\eta$ that is concentrated on $\csg(\Sp)$.
\end{theorem}

\begin{proof}
	Fix $\eta\in \P(\CC)$ with $(e_0)_\#\eta = \mu$ and $(e_1)_\#\eta =\nu$. Then $(e_0\times e_1)_\#\eta\in \Gamma(\mu,\nu)$ and, by Lemma \ref{lem_Ainfty}, we have
	\begin{align}\label{eq_Aineq}
		\norm{\A}_{L^\infty(\eta)} \geq \norm{\d\circ(e_0\times e_1)}_{L^\infty(\eta)} = \norm{\d(\cdot,\cdot)}_{L^\infty((e_0\times e_1)_\#\eta)} \geq W_\infty(\mu,\nu).
	\end{align}
	This proves one inequality and the characterization of minimizers. To prove the other, let $\pi\in \Gamma_\infty(\mu,\nu)$ be optimal and let $\Phi:\Sp\times\Sp\to \csg(\Sp)$ be a $\pi$-measurable map such that $\omega_{x,y}:=\Phi(x,y)\in \csg(\Sp)$ satisfies $\omega_{x,y}(0)=x$ and $\omega_{x,y}(1)=y$, and define $\eta:=\Phi_\#\pi$. Then, clearly $\eta$ is concentrated on the set of constant speed geodesics. Therefore, with this choice of $\eta$, all the inequalities in \eqref{eq_Aineq} are actually equalities.
	The characterization of minimizers now easily follows.

\end{proof}

\begin{remark}
	In the $q$-Wasserstein space, with $1\le q\le \infty$, property \ref{th_DynForm_i} of Theorem \ref{th_DynamicFormulation} is replaced by (the equivalent condition) \emph{“$\eta$ is concentrated on $\csg(\Sp)$”}. However, in the $\infty$-Wasserstein space there may exist minimizers $\eta$ for the dynamic transportation problem \eqref{th_DynForm_claim} that are not concentrated on the family of constant speed geodesics. As an example, consider the space $\Sp= \R^2$ and the measures 
	\begin{gather*}
	\mu:=\frac{1}{4}\L^2\llcorner [-1/2,3/2]\times[-1,1],\\
	 \nu:=\frac{1}{4}\left(\L^2\llcorner[-1/2,1/2]\times [-1,1] + \L^2\llcorner [9/2,5]\times[-1,1]\right).
	\end{gather*}
	Let $\eta\in \P(\C^0([0,1];\R^2))$ be the measure defined by duality with $\C^0_b(\C^0([0,1];\Rn))$ as
	\begin{align*}
	\int_{\C^0([0,1];\Rn)}\psi(\omega)\,d\eta(\omega) :=\frac{1}{4}\int_{-\frac{1}{2}}^\frac{1}{2}\int_{-1}^1 \psi(R_{x,y})\,dy\,dx + \frac{1}{4}\int_{\frac{1}{2}}^\frac{3}{2}\int_{-1}^1 \psi(T_{x,y})\,dy\,dx,
	\end{align*}
	where for every $(x,y)\in \R^2$, $R_{x,y},T_{x,y}:[0,1]\to \R^2$ and  are the paths
	\begin{gather*}
	R_{x,y}(t) := (x\,\cos(\pi t) - y\sin(\pi t), x\,\sin(\pi t) + y\cos(\pi t) )\quad \forall 0\leq t\leq 1,\\
	T_{x,y}(t) := (x + 4t, y)\quad \forall 0\leq t\leq 1.
	\end{gather*}
	Clearly $(e_0)_\#\eta = \mu$ and $(e_1)_\#\eta =\nu$. Moreover, as
	\begin{gather*}
	|\dot R_{x,y}(t)| = \pi\sqrt{x^2 + y^2} \leq \frac{\pi \sqrt{5}}{2} \le 4\quad \forall (x,y)\in [-1/2,1/2]\times[-1,1],\\
	|\dot T_{x,y}(t)| = 4 \quad \forall (x,y)\in [1/2,3/2]\times[-1,1],
	\end{gather*}
	it follows that $\norm{\A}_{L^\infty(\eta)}= 4$. On the other hand it is easy to see (use for instance \eqref{eq_Winfty}) that $W_\infty(\mu,\nu)=4$. Therefore $\eta$ is a minimizer for the dynamic transportation problem that gives positive measure to the set $\{R_{x,y}:-1/2\leq x\leq 1/2,\,-1\leq y\leq 1\}\minus\{R_{0,0}\}$, and the latter has empty intersection with $\csg(\R^2)$. 
\end{remark}

\end{appendices}

\renewcommand{\refname}{References}
\markboth{references}{references}

{\bbsizz{

}}


\begin{thebibliography}{10}

\bibitem{AldrigoInPreparation}
{\sc Aldrigo, P.},
\newblock Connection by AC arcs of the space PL$_\infty^p$,
\newblock in preparation.

\bibitem{Ambrosio1995}
{\sc Ambrosio, L.}
\newblock Minimizing movements.
\newblock {\em Rendiconti dell'Accademia Nazionale delle Scienze detta dei XL.
  Memorie di Matematica e Applicazioni 19\/} (1995), 191--246.

\bibitem{ambrosio2004lecture}
{\sc Ambrosio, L.}
\newblock {Lecture notes on optimal transport problems}.
\newblock In {\em Mathematical Aspects of Evolving Interfaces: Lectures given
  at the CIM-CIME joint Euro-Summer School held in Madeira, Funchal, Portugal,
  July 3-9, 2000}. Springer, 2004, pp.~1--52.

\bibitem{ambrosio2021}
{\sc Ambrosio, L., Bru{\'e}, E., Semola, D., et~al.}
\newblock {\em Lectures on optimal transport}, vol.~130.
\newblock Springer, 2021.

\bibitem{ambrosio2005gradient}
{\sc Ambrosio, L., Gigli, N., and Savar{\'e}, G.}
\newblock {\em {Gradient Flows: in Metric Spaces and in the Space of
  Probability Measures}}.
\newblock Springer, 2005.

\bibitem{AmbrosioKirchheim2000}
{\sc Ambrosio, L., and Kirchheim, B.}
\newblock {Rectifiable Sets in Metric and Banach Spaces}.
\newblock {\em Mathematische Annalen 318}, 3 (2000), 527--555.

\bibitem{Aubin1976}
{\sc Aubin, T.}
\newblock {Probl{\`e}mes isop{\'e}rim{\'e}triques et espaces de Sobolev}.
\newblock {\em Journal of Differential Geometry 11}, 4 (1976), 573--598.

\bibitem{BenamouBrenier2000}
{\sc Benamou , J.D., and Brenier, Y.}
\newblock {A computational fluid mechanics solution to the Monge--Kantorovich mass transfer problem},
\newblock {\em Numerische Mathematik  84}, (2000), no.~3, 375--393.

\bibitem{Brasco2010}
{\sc Brasco, L.}
\newblock A survey on dynamical transport distances.
\newblock {\em Journal of Mathematical Sciences\/} (2010).

\bibitem{CailletSantambrogio2024}
{\sc Caillet, T., and Santambrogio, F.}
\newblock {Doubly Nonlinear Diffusive {PDEs}: New Existence Results via
  Generalized Wasserstein Gradient Flows}.
\newblock {\em SIAM Journal on Mathematical Analysis 56}, 6 (2024), 7043--7073.

\bibitem{CarlierPoon2019}
{\sc Carlier, G., and Poon, C.}
\newblock {On the total variation Wasserstein gradient flow and the {TV-JKO}
  scheme}.
\newblock {\em ESAIM: Control, Optimisation and Calculus of Variations 25\/}
  (2019), Paper No. 42.

\bibitem{ChambolleDuvalMachado2023}
{\sc Chambolle, A., Duval, V., and Machado, J.~M.}
\newblock {The Total Variation--Wasserstein Problem: A New Derivation of the
  Euler--Lagrange Equations}.
\newblock In {\em Geometric Science of Information\/} (2023), vol.~14071 of
  {\em Lecture Notes in Computer Science}, Springer Nature Switzerland,
  pp.~610--619.

\bibitem{champion2008wasserstein}
{\sc Champion, T., De~Pascale, L., and Juutinen, P.}
\newblock {The $\infty$-Wasserstein distance: local solutions and existence of
  optimal transport maps}.
\newblock {\em SIAM Journal on Mathematical Analysis 40}, 1 (2008), 1--20.

\bibitem{CorderoErausquinNazaretVillani2004}
{\sc Cordero-Erausquin, D., Nazaret, B., and Villani, C.}
\newblock {A mass-transportation approach to sharp Sobolev and
  Gagliardo--Nirenberg inequalities}.
\newblock {\em Advances in Mathematics 182}, 2 (2004), 307--332.

\bibitem{DeGiorgi1958}
{\sc De~Giorgi, E.}
\newblock Sulla propriet{\`a} isoperimetrica dell'ipersfera, nella classe degli
  insiemi aventi frontiera orientata di misura finita.
\newblock {\em Atti della Accademia Nazionale dei Lincei. Memorie. Classe di
  Scienze Fisiche, Matematiche e Naturali. Sezione I 5\/} (1958), 33--44.

\bibitem{DeGiorgi1993}
{\sc De~Giorgi, E.}
\newblock New problems on minimizing movements.
\newblock In {\em Boundary Value Problems for PDE and Applications}. Masson,
  1993, pp.~81--98.

\bibitem{EvansGariepy2015}
{\sc Evans, L.~C., and Gariepy, R.~F.}
\newblock {\em Measure Theory and Fine Properties of Functions}, revised~ed.
\newblock CRC Press, 2015.

\bibitem{figalli2021invitation}
{\sc Figalli, A., and Glaudo, F.}
\newblock {\em An invitation to optimal transport, Wasserstein distances, and
  gradient flows}, vol.~1.
\newblock EMS press Berlin, 2021.

\bibitem{FigalliMaggiPratelli2010}
{\sc Figalli, A., Maggi, F., and Pratelli, A.}
\newblock A mass transportation approach to quantitative isoperimetric
  inequalities.
\newblock {\em Inventiones Mathematicae 182}, 1 (2010), 167--211.

\bibitem{FigalliMaggiPratelli2013}
{\sc Figalli, A., Maggi, F., and Pratelli, A.}
\newblock {Sharp stability theorems for the anisotropic Sobolev and log-Sobolev
  inequalities on functions of bounded variation}.
\newblock {\em Advances in Mathematics 242\/} (2013), 80--101.

\bibitem{FlemingRishel1960}
{\sc Fleming, W.~H., and Rishel, R.}
\newblock An integral formula for total gradient variation.
\newblock {\em Archiv der Mathematik 11\/} (1960), 218--222.

\bibitem{Galdi2011}
{\sc Galdi, G.~P.}
\newblock {\em An Introduction to the Mathematical Theory of the Navier--Stokes
  Equations}, second~ed.
\newblock Springer Monographs in Mathematics. Springer, New York, 2011.

\bibitem{givens}
{\sc Givens, C.~R., and Shortt, R.~M.}
\newblock {A class of Wasserstein metrics for probability distributions.}
\newblock {\em Michigan Mathematical Journal 31}, 2 (1984), 231--240.

\bibitem{JordanKinderlehrerOtto1998}
{\sc Jordan, R., Kinderlehrer, D., and Otto, F.}
\newblock The variational formulation of the {Fokker--Planck} equation.
\newblock {\em SIAM Journal on Mathematical Analysis 29}, 1 (1998), 1--17.

\bibitem{LinSantambrogio2026}
{\sc Lin, K., and Santambrogio, F.}
\newblock Existence of a solution of the {TV} {Wasserstein} gradient flow,
  2026, ArXiv: https://doi.org/10.48550/arXiv.2601.22847.

\bibitem{Lisini2007}
{\sc Lisini, S.}
\newblock {Characterization of absolutely continuous curves in Wasserstein
  spaces}.
\newblock {\em Calculus of Variations and Partial Differential Equations 28\/}
  (2007), 85--120.

\bibitem{santambrogio20151}
{\sc Santambrogio, F.}
\newblock {\em Optimal Transport for Applied Mathematicians: Calculus of
  Variations, PDEs, and Modeling}, vol.~1.
\newblock Springer, 2015.

\bibitem{Santambrogio2017}
{\sc Santambrogio, F.}
\newblock {Euclidean, metric, and Wasserstein gradient flows: an overview}.
\newblock {\em Bulletin of Mathematical Sciences 7}, 1 (2017), 87--154.

\bibitem{Talenti1976}
{\sc Talenti, G.}
\newblock {Best constant in Sobolev inequality}.
\newblock {\em Annali di Matematica Pura ed Applicata 110\/} (1976), 353--372.

\end{thebibliography}
\end{document}